\documentclass[a4paper]{article}
\usepackage{amssymb}
\usepackage{amsthm}
\usepackage{tikz}
\usepackage{mathtools}
\usepackage[title,toc,page,header]{appendix}
\usepackage{titletoc}
\usepackage{dutchcal}

\usepackage[mathscr]{euscript} 
\DeclareFontFamily{OT1}{pzc}{}
\DeclareFontShape{OT1}{pzc}{m}{it}{<-> s * [1.10] pzcmi7t}{}
\DeclareMathAlphabet{\mathpzc}{OT1}{pzc}{m}{it}

\dottedcontents{section}[2.5em]{\bfseries}{2.5em}{1pc}

\usetikzlibrary{calc}
\usetikzlibrary{matrix,arrows,decorations}
\newtheorem{theorem}{Theorem}[section]
\newtheorem{lemma}[theorem]{Lemma}

\newtheorem{remark}[theorem]{Remark}
\newtheorem{definition}[theorem]{Definition} 
\newtheorem{example}[theorem]{Example}
\newtheorem{construction}[theorem]{Construction}

\usepackage{hyperref}

\title{Categories of vector spaces and Grassmannians}

\author{Yi-Sheng Wang}
\date{}

\begin{document}
\maketitle
\begin{abstract}
We construct a new category of vector spaces which contains both the standard category of vector spaces and Grassmannians. Its space of objects classifies vector bundles, its space of morphisms classifies bundle isomorphisms, and it can be delooped to obtain the topological $K$-theory spectrum. 
\end{abstract}
\tableofcontents
\addtocontents{toc}{\protect\setcounter{tocdepth}{1}}
\section{Introduction}
There are many different models for the classifying space of the general linear group $\operatorname{GL}_{k}(\mathbb{F})$, where $\mathbb{F}=\mathbb{R}$ or $\mathbb{C}$. Two most commonly used models are the Grassmannian manifold $\operatorname{Gr}_{k}(\mathbb{F}^{\infty})$ and the geometric realization of the bar construction $B_{\cdot}\operatorname{GL}_{k}(\mathbb{F})$, denoted by $B\operatorname{GL}_{k}(\mathbb{F})$ \cite{BV},\cite{Ma1}, \cite{Ma3}. These two models are extremely important in the study of algebraic and geometric topology. The Grassmannian $\operatorname{Gr}_{k}(\mathbb{F}^{\infty})$ is more concrete and geometric; it is a colimit of a sequence of compact smooth manifolds, whose combinatorial and differential structures are well understood. On the other hand, the space $B\operatorname{GL}_{k}(\mathbb{F})$, though more abstract, fits better into the framework of delooping machines and hence is useful in the infinite loop space theory. As both spaces classify the isomorphism classes of $k$-dimensional vector bundles over paracompact Hausdorff spaces, these two spaces have the same homotopy type. However, there is no concrete description of the homotopy equivalence between them in the literature, to the author's knowledge. In this paper, we construct a new category of vector spaces, which contains both the geometric features of the Grassmannian $\operatorname{Gr}_{k}(\mathbb{F}^{\infty})$ and the algebraic properties of the space $B\operatorname{GL}_{k}(\mathbb{F})$. Using this category, we obtain a concrete description of a zig-zag of homotopy equivalences between $\operatorname{Gr}_{k}(\mathbb{F}^{\infty})$ and $B\operatorname{GL}_{k}(\mathbb{F})$.

Throughout the paper, we work in the category of weak Hausdorff $k$-spaces, denoted by $\mathpzc{Top}^{w}$.  
\subsection*{Statements of the results} 
Let $\mathcal{G}$ be the topological groupoid in $\mathpzc{Top}^{w}$ whose spaces of objects and morphisms are the disjoint union of the Grassmannian manifolds 
\[\coprod_{k}\operatorname{Gr}_{k}(\mathbb{F}^{\infty}),\] 
and let $\mathcal{V}_{\mathbb{F}}$ be the enriched category in $\mathpzc{Top}^{w}$ whose set of objects and the space of morphisms are the discrete set $\{\mathbb{F}^{n}\}_{n\in\mathbb{N}\cup\{0\}}$ and the disjoint union of the spaces of $m$-by-$n$ matrices $\coprod\limits_{m,n}M_{m,n}(\mathbb{F})$, respectively. Now, consider the vector bundle  
\[(s^{m},t^{n}):\coprod_{k,l}\operatorname{Mor}_{k,l}^{m,n}(\mathbb{F})\rightarrow \coprod_{k}\operatorname{Gr}_{k}(\mathbb{F}^{m})\times\coprod_{l}\operatorname{Gr}_{l}(\mathbb{F}^{n})\]  
whose total space consists of linear maps from $V$ to $W$, for any $V\in\coprod\limits_{k}\operatorname{Gr}_{k}(\mathbb{F}^{m})$ and $W\in\coprod\limits_{l}\operatorname{Gr}_{l}(\mathbb{F}^{n})$ (see \eqref{CoordinateMor} for its coordinate charts), and observe that the identity of a vector space and the composition of linear maps give an inclusion  
\[e^{m}:\coprod_{k}\operatorname{Gr}_{k}(\mathbb{F}^{m})\rightarrow \coprod_{k,k}\operatorname{Mor}_{k,k}^{m,m}(\mathbb{F})\] 
and a natural map
\[\circ^{m,n,p}: \coprod_{k,l,j}\operatorname{Mor}_{k,l}^{m,n}(\mathbb{F})\times_{\operatorname{Gr}_{l}(\mathbb{F}^{n})}\operatorname{Mor}_{l,j}^{n,p}(\mathbb{F})\rightarrow \coprod_{k,j}\operatorname{Mor}_{k,j}^{m.p}(\mathbb{F}),\]
respectively. Then we have the following theorem.
\begin{theorem}[Lemma \ref{HsemiringVFf}; Lemma \ref{Manyrealizations}; Theorem \ref{ComparisonVFVFfG}]\label{MainThm}
Define
\begin{align*}
O(\mathcal{V}_{\mathbb{F}}^{f})&:=\coprod_{k}\operatorname{Gr}_{k}(\mathbb{F}^{\infty});\\
M(\mathcal{V}_{\mathbb{F}}^{f})&:=\coprod_{k,l}\operatorname*{colim}\limits_{m,n}\operatorname{Mor}_{k,l}^{m,n}(\mathbb{F});\\
s&:=\operatorname*{colim}\limits_{m}s^{m};\\
t&:=\operatorname*{colim}\limits_{m}t^{m};\\
e&:=\operatorname*{colim}\limits_{m}e^{m};\\
\circ&:=\operatorname*{colim}\limits_{m}\circ^{m,m,m}.
\end{align*}
Then the $6$-tuple 
\[\mathcal{V}_{\mathbb{F}}^{f}:=\{O(\mathcal{V}_{\mathbb{F}}^{f}),M(\mathcal{V}_{\mathbb{F}}^{f}), s, t, e, \circ\}\]
constitutes an internal category in $\mathpzc{Top}^{w}$, and the canonical inclusions
\[\mathcal{V}_{\mathbb{F}}\rightarrow \mathcal{V}_{\mathbb{F}}^{f}\leftarrow \mathcal{G}\] 
induce a cospan of homotopy equivalences of $\operatorname{H}$-semiring spaces
\[\vert\operatorname{Ner}_{\cdot}iso\mathcal{V}_{\mathbb{F}}\vert\rightarrow \vert\operatorname{Ner}_{\cdot}iso\mathcal{V}_{\mathbb{F}}^{f}\vert\leftarrow \vert\operatorname{Ner}_{\cdot}\mathcal{G}\vert,\]
where $iso\mathcal{C}$ is the internal subcategory of isomorphisms in $\mathcal{C}$, $\operatorname{Ner}_{\cdot}$ is the nerve construction and $\vert-\vert$ is the geometric realization functor. Furthermore, the singularization of the category $\mathcal{V}_{\mathbb{F}}^{f}$ 
is a simplicial $S$-category equipped with a sum functor \cite[Definition $1.18$]{Wang5}, and its associated $\Omega$-prespectrum represents topological $K$-theory.
\end{theorem}

\subsection*{Outline of the paper}
In Section $2$, following Boardman and Vogt's approach \cite{BV} \cite[p.210-213]{Sw}, we give a detailed account of the $H$-semiring structure on the disjoint union of Grassmannians
\[\coprod\limits_{k}\operatorname{Gr}_{k}(\mathbb{F}^{\infty})\]
and its topological group completion. Section $3$ reviews some properties of the category $\mathcal{V}_{\mathbb{F}}$. In Section $4$, we examine carefully the construction of the category $\mathcal{V}_{\mathbb{F}}^{f}$, including its internal structures, the geometric meaning of these structures and its nerve. The last section discusses a comparison between these three categories $\mathcal{G}$, $\mathcal{V}_{\mathbb{F}}$ and $\mathcal{V}_{\mathbb{F}}^{f}$. In particular, the proof of Theorem \ref{MainThm} is completed in that section. In appendix $A$, we collect some properties of weak Hausdorff spaces needed for this paper, including the commutativity of pullbacks and colimits, and various constructions that preserve cofibrations. In Appendix $B$, some basic concepts of internal category theory are reviewed. 
  
\subsection*{Notation and conventions}
Throughout the paper, all topological spaces are weak Hausdorff $k$-spaces, unless otherwise specified. Given two topological spaces $X$ and $Y$, the set of homotopy classes of maps from $X$ to $Y$ is denoted by $[X,Y]$. When both $X$ and $Y$ are pointed, we denote the set of homotopy classes of pointed maps by $[X,Y]_{\ast}$. Given an internal category $\mathcal{C}$ in $\mathpzc{Top}^{w}$, $O(\mathcal{C})$ and $M(\mathcal{C})$ denote its spaces of objects and morphisms, respectively.

\section{Grassmannians}
\subsection{$H$-semiring structure on $\coprod\limits_{k}\operatorname{Gr}_{k}(\mathbb{F}^{\infty})$}
This subsection examines the $H$-semiring structure on $\coprod\limits_{k}\operatorname{Gr}_{k}(\mathbb{F}^{\infty})$, where $\operatorname{Gr}_{k}(\mathbb{F}^{\infty})$ is the Grassmannian manifold of $k$-dimensional vector subspaces in $\mathbb{F}^{\infty}$. We first recall a useful lemma due to Boardman and Vogt \cite{BV}:
\begin{lemma}\label{BVlemma}
Given an inner product vector space $V$ over $\mathbb{F}$ of countable dimension, then the space of linear isomerties $I(V,\mathbb{F}^{\infty})$ is contractible, where $\mathbb{F}^{\infty}$ admits a natural inner product structure. 
\end{lemma}
\begin{proof}
See \cite[p.212]{Sw}
\end{proof}
\noindent
In particular, Lemma \ref{BVlemma} implies that any two isometries 
\[f,g:V\rightarrow\mathbb{F}^{\infty}\]
are homotopic through isometries. With this lemma, we give a detailed proof of the following well-known result.
\begin{theorem}\label{HsemiringstronG}
The tensor product and direct sum of vector spaces induce a $H$-semiring structure on $\coprod\limits_{k}\operatorname{Gr}_{k}(\mathbb{F}^{\infty})$.
\end{theorem}
\begin{proof}  
To construct the additive structure, we choose an isometric isomorphism 
\[\alpha:\mathbb{F}^{\infty}\oplus\mathbb{F}^{\infty}\rightarrow\mathbb{F}^{\infty},\]
and define the additive operation $a$ to be the composition:
\[a:\operatorname{Gr}_{k}(\mathbb{F}^{\infty})\times \operatorname{Gr}_{l}(\mathbb{F}^{\infty})\rightarrow\operatorname{Gr}_{k+l}(\mathbb{F}^{\infty}\oplus\mathbb{F}^{\infty})\xrightarrow{\operatorname{Gr}_{k+l}(\alpha)}\operatorname{Gr}_{k+l}(\mathbb{F}^{\infty}).\]
The associativity and commutativity of the operation $a$ follows from Lemma \ref{BVlemma} and the homotopy commutative diagrams below---for the sake of simplicity, we use $\mathbb{F}$ to denote $\operatorname{Gr}_{k}(\mathbb{F}^{\infty})$ and $\mathbb{F}\oplus \mathbb{F}$ to denote $\operatorname{Gr}_{k}(\mathbb{F}^{\infty}\oplus\mathbb{F}^{\infty})$; the indicies $k$ and $l$ are also omitted as it is easy to fill them in.

\noindent
\textbf{Associativity:} 
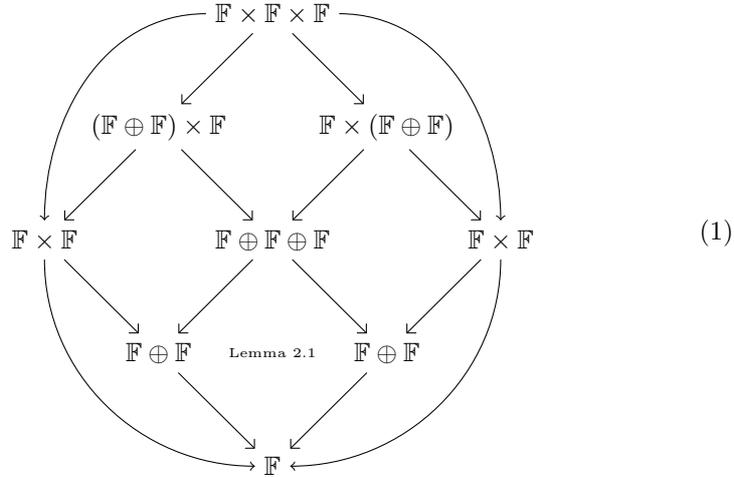
\begin{center}
\begin{equation}\label{AssoDia} 
\begin{tikzpicture}[baseline=(current  bounding  box.center)]

\node (U) at (3,3) {$\mathbb{F}\times \mathbb{F}\times \mathbb{F}$};
\node (UMl) at (1.5,1.5) {$(\mathbb{F}\oplus\mathbb{F})\times\mathbb{F}$}; 
\node (UMr) at (4.5,1.5) {$\mathbb{F}\times(\mathbb{F}\oplus\mathbb{F})$};
\node (Ml) at (0,0) {$\mathbb{F}\times\mathbb{F}$};
\node (Mm) at (3,0) {$\mathbb{F}\oplus\mathbb{F}\oplus\mathbb{F}$};
\node (Mr) at (6,0) {$\mathbb{F}\times\mathbb{F}$};
\node (LMl) at (1.5,-1.5) {$\mathbb{F}\oplus\mathbb{F}$};
\node (LMr) at (4.5,-1.5) {$\mathbb{F}\oplus\mathbb{F}$};
\node (L) at (3,-3){$\mathbb{F}$};

\node (Label) at (3,-1.5){\tiny Lemma \ref{BVlemma}};

\path[->, font=\scriptsize,>=angle 90] 

(U) edge (UMl)
(U) edge (UMr)
(UMl) edge (Ml)
(UMl) edge (Mm)
(UMr) edge (Mm)
(UMr) edge (Mr)
(Ml) edge (LMl)
(Mm) edge (LMl)
(Mm) edge (LMr)
(Mr) edge (LMr)
(LMl) edge (L)
(LMr) edge (L);

\draw [->] (U) to [out=0, in=90] (Mr);
\draw [->] (U) to [out=180,in=90] (Ml);
\draw [->] (Mr) to [out=-90,in=0] (L);
\draw [->] (Ml) to [out=-90,in=180] (L);
\end{tikzpicture} 
\end{equation}
\end{center}
\textbf{Commutativity:}
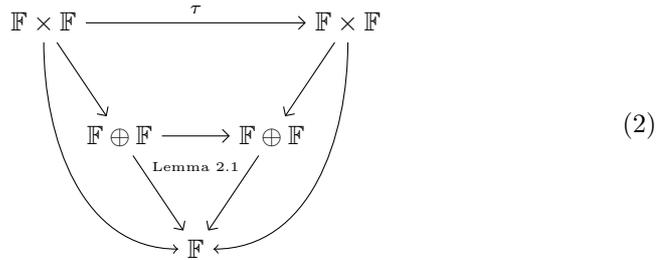
\begin{center}
\begin{equation}\label{CommDia}
\begin{tikzpicture}[baseline=(current  bounding  box.center)]
\node (Lu) at (0,3){$\mathbb{F}\times\mathbb{F}$};
\node (Ru) at (4,3){$\mathbb{F}\times\mathbb{F}$};
\node (Lm) at (1,1.5){$\mathbb{F}\oplus\mathbb{F}$};
\node (Rm) at (3,1.5){$\mathbb{F}\oplus\mathbb{F}$};
\node (Ml) at (2,0){$\mathbb{F}$};
\node (Label) at (2,1.1){\tiny Lemma \ref{BVlemma}};

\path[->, font=\scriptsize,>=angle 90] 

(Lu) edge (Lm)
(Lu) edge node [above]{$\tau$}(Ru)
(Ru) edge (Rm)
(Lm) edge (Rm)
(Lm) edge (Ml)
(Rm) edge (Ml);

\draw [->] (Lu) to [out=-90,in=180] (Ml);
\draw [->] (Ru) to [out=-90,in=0] (Ml);

\end{tikzpicture}
\end{equation}
\end{center}
where $\tau$ is given by the assignment $(x,y)\mapsto (y,x)$. With the same notation, the following homotopy commutative diagram  
\begin{center}
\begin{equation}\label{AddUnit}
\begin{tikzpicture}[baseline=(current  bounding  box.center)] 
\node (Lu) at (0,3) {$\mathbb{F}$};
\node (Ru) at (3,3) {$\mathbb{F}\times \mathbb{F}$};
\node (Rm) at (3,1.5) {$\mathbb{F}\oplus\mathbb{F}$};
\node (Rl) at (3,0) {$\mathbb{F}$};

\node (Label) at (1.5,1){\tiny Lemma \ref{BVlemma}};

\path[->, font=\scriptsize,>=angle 90] 
(Lu) edge node [above] {$i_{1}/i_{2}$}(Ru)
(Ru) edge (Rm)
(Rm) edge (Rl)
(Lu) edge (Rm);

\draw [->] (Lu) to [out=-90,in=180] (Rl);
\draw [->] (Ru) to [out=-45,in=45] (Rl);

\end{tikzpicture}
\end{equation}
\end{center}
shows $\ast=\operatorname{Gr}_{0}(\mathbb{F}^{\infty})$ is an additive identity, where $\iota_{1}(x)=(x,\ast)$ and $\iota_{2}(x)=(\ast,x)$.

Similarly, one can define the multiplicative operation $m$ to be the composition: 
\[m: \operatorname{Gr}_{k}(\mathbb{F}^{\infty})\times\operatorname{Gr}_{l}(\mathbb{F}^{\infty})\rightarrow \operatorname{Gr}_{k+l}(\mathbb{F}^{\infty}\otimes\mathbb{F}^{\infty})\xrightarrow{\operatorname{Gr}_{k+l}(\beta)}\operatorname{Gr}_{k+l}(\mathbb{F}^{\infty}),\]
where $\beta$ is an isometric isomorphism $\mathbb{F}^{\infty}\otimes\mathbb{F}^{\infty}\rightarrow \mathbb{F}^{\infty}$. The associativity and commutativity of $m$ follow from the diagrams \eqref{AssoDia} and \eqref{CommDia} with $\oplus$ replaced by $\otimes$ and $\alpha$ by $\beta$; letting $i_{1}(x)=x\otimes e_{1}$ and $i_{2}(x)=e_{1}\otimes x$, then the diagram \eqref{AddUnit} with $\oplus$ replaced by $\otimes$ implies that the point  
\[ \{\mathbb{F}\}=\operatorname{Gr}_{1}(\mathbb{F})\hookrightarrow\operatorname{Gr}_{1}(\mathbb{F}^{\infty}).\]
is a multiplicative identity, where $\{e_{i}\}_{i\in\mathbb{N}}$ is the standard basis of $\mathbb{F}^{\infty}$. 
 
Lastly, the distributivity of $m$ with respect to $a$ results from the homotopy commutative diagram  
\begin{center} 
\begin{tikzpicture}
\node (LL) at (4.5,0) {$\mathbb{F}$};
\node (MLl) at (1.5,1.5) {$\mathbb{F}\oplus\mathbb{F}$};
\node (MLr) at (7.5,1.5) {$\mathbb{F}\otimes\mathbb{F}$};
\node (Ml) at (0,3) {$\mathbb{F}\times\mathbb{F}$};
\node (Mm1)at (3,3) {$(\mathbb{F}\otimes\mathbb{F})\oplus(\mathbb{F}\otimes\mathbb{F})$};
\node (Mm2) at (6,3) {$(\mathbb{F}\oplus\mathbb{F})\otimes\mathbb{F}$};
\node (Mr) at (9,3) {$\mathbb{F}\times\mathbb{F}$};
\node (MUl) at (1.5,4.5){$\mathbb{F}\otimes\mathbb{F}\times\mathbb{F}\otimes\mathbb{F}$};
\node (MUr) at (7.5,4.5){$\mathbb{F}\oplus\mathbb{F}\times\mathbb{F}$};
\node (Extra) at (3,6){$\mathbb{F}\times\mathbb{F}\times\mathbb{F}\times\mathbb{F}$};
\node (UU) at (4.5,7.5) {$\mathbb{F}\times\mathbb{F}\times\mathbb{F}$};

\node (Label) at (4.5,1.5){\tiny Lemma \ref{BVlemma}};

\path[->, font=\scriptsize,>=angle 90] 

(UU) edge (Extra)
(UU) edge (MUr)
(Extra) edge (MUl)
(MUr) edge (Mm2)
(MUl) edge (Mm1)
(MUr) edge (Mr)
(MUl) edge (Ml)
(Mr) edge (MLr)
(Ml) edge (MLl)
(Mm2) edge (MLr)
(Mm1) edge (MLl)
  
(MLl) edge (LL)
(MLr) edge (LL);

\draw [double equal sign distance] (Mm1) to (Mm2);
\draw [->] (UU) to [out=0,in=90] (Mr);
\draw [->] (UU) to [out=180,in=90] (Ml);

\draw [->] (Ml) to [out=-90,in=180] (LL);
\draw [->] (Mr) to [out=-90,in=0] (LL);
\end{tikzpicture}
\end{center}  
The right distribution law can be proved similarly.

\end{proof}
\begin{remark} 
One can choose the homotopies in the homotopy commutative diagrams for multiplicative operation $b$ and multiplicative identity $e_{1}$ in such a way that the homotopies preserve the multiplicative identity. 
\end{remark}

\begin{remark}
It is well-known that $\coprod\limits_{k}\operatorname{Gr}_{k}(\mathbb{F}^{\infty})$ is homotopy equivalent to an $E_{\infty}$-ring space \cite[Example $5.4$ and Section VII.$4$]{Ma4}. Instead of studying its higher structure, our goal here is to give a concrete and elementary description of its $H$-semiring space structure. 
\end{remark}

\subsection{The topological group completion of $\coprod\limits_{k}\operatorname{Gr}_{k}(\mathbb{F}^{\infty})$}\label{GroupCompletionthesection}
In this section, we give a detailed construction of the topological group completion of $\coprod\limits_{k}\operatorname{Gr}_{k}(\mathbb{F}^{\infty})$. 

Observe first that the assignment $V\mapsto \alpha_{\ast}V+e_{1}$ induces a map $\iota^{\prime}_{k}:\operatorname{Gr}_{k}(\mathbb{F}^{\infty})\rightarrow \operatorname{Gr}_{k+1}(\mathbb{F}^{\infty})$, where $\alpha$ is given by
\begin{align*}
\alpha: \mathbb{F}^{\infty}&\rightarrow \mathbb{F}^{\infty}\\  
(x_{1},...,x_{n},...)&\mapsto (0,x_{1},...,x_{n},...),
\end{align*}
and it is well-known that the colimit of the sequence $\{\coprod\limits_{k}\operatorname{Gr}_{k}(\mathbb{F}^{\infty}),\iota^{\prime}_{k}\}$
has the correct homotopy type of the topological group completion of $\coprod\limits_{k}\operatorname{Gr}_{k}(\mathbb{F}^{\infty})$. However, it is not easy to describe the $H$-space structure on $\operatorname*{colim}\limits_{k}\operatorname{Gr}_{k}(\mathbb{F}^{\infty})$. Thus, we employ another construction introduced in \cite[p.210-3]{Sw}. Denote the Grassmannian of $k$-dimensional vector subspaces in $\mathbb{F}^{k}\otimes\mathbb{F}^{\infty}$ by 
\[G(\mathbb{F}^{k}):=\operatorname{Gr}_{k}(\mathbb{F}^{k}\otimes\mathbb{F}^{\infty}),\]
and observe that, for each $k$, there is a natural cofibration given by 
\begin{align*}
\iota_{k}:G(\mathbb{F}^{k})&\hookrightarrow G(\mathbb{F}^{k+1}),\\
 V  &\mapsto V+<e_{k+1}\otimes e_{1}>,
\end{align*}
where $\{e_{i}\}$ is the standard basis of $\mathbb{F}^{k}$ (resp. $\mathbb{F}^{\infty}$). Let $G(\mathbb{F}^{\infty})$ to be the colimit of the sequence $\{G(\mathbb{F}^{k}),\iota_{k}\}$. Then there is a natural map  
\[\coprod_{k}G(\mathbb{F}^{k})\rightarrow \mathbb{Z}\times G(\mathbb{F}^{\infty}).\]

\begin{lemma}\label{Ginfinity}
Both $\coprod\limits_{k}G(\mathbb{F}^{k})$ and $\mathbb{Z}\times G(\mathbb{F}^{\infty})$
are $H$-spaces and the canonical map
\[\coprod_{k}G(\mathbb{F}^{k})\rightarrow \mathbb{Z}\times G(\mathbb{F}^{\infty})\]
is a $H$-map.
\end{lemma}
\begin{proof}
The approaches presented here is taken from \cite[p.211-213]{Sw}. The construction of the $H$-space structure on $\coprod\limits_{k}G(\mathbb{F}^{k})$ is very similar to that on $\coprod\limits_{k}\operatorname{Gr}_{k}(\mathbb{F}^{\infty})$. The additive operation 
\[a:\coprod_{k}G(\mathbb{F}^{k})\times\coprod_{k}G(\mathbb{F}^{k})\rightarrow \coprod_{k}G(\mathbb{F}^{k})\]  
is induced by the composition
\[G(\mathbb{F}^{k})\times G(\mathbb{F}^{l})\rightarrow G(\mathbb{F}^{k}\oplus\mathbb{F}^{l})\xrightarrow{G(\alpha)} G(\mathbb{F}^{k+l}),\]
where $\alpha$ is an isometric isomorphism
\[\mathbb{F}^{k}\oplus\mathbb{F}^{l}\xrightarrow{\alpha}\mathbb{F}^{k+l}.\] 
It is clear that the point $G(\mathbb{F}^{0})=\ast$ is the additive identity with respect to the additive operation by the following homotopy commutative diagram, where $\mathbb{F}^{k}$ stands for $G(\mathbb{F}^{k})$: 
\begin{center}
\begin{tikzpicture}

\node (Lu) at (0,3) {$\mathbb{F}^{k}$};
\node (Ru) at (3,3) {$\mathbb{F}^{k}\times\mathbb{F}^{0}/\mathbb{F}^{0}\times \mathbb{F}^{k}$};
\node (Rm) at (3,1.5) {$\mathbb{F}^{k}\oplus\mathbb{F}^{0}/\mathbb{F}^{0}\oplus\mathbb{F}^{k}$};
\node (Rl) at (3,0) {$\mathbb{F}^{k}$};

\node (Label) at (1.5,1){\tiny Lemma \ref{BVlemma}};

\path[->, font=\scriptsize,>=angle 90] 
(Lu) edge node [above] {$i_{1}/i_{2}$}(Ru)
(Ru) edge (Rm)
(Rm) edge (Rl)
(Lu) edge (Rm);

\draw [->] (Lu) to [out=-90,in=180] (Rl);
\draw [->](Ru.350) to [out=-45,in=45] (Rl);

\end{tikzpicture}
\end{center}
On the other hand, the following homotopy commutative diagrams, which are essentially the diagrams \eqref{AssoDia} and \eqref{CommDia}, show the commutativity and associativity of the additive operation $a$
\begin{center}
\begin{equation}\label{AssoDia2}
\begin{tikzpicture}[baseline=(current bounding box.center)]
\node (U) at (3,3) {$\mathbb{F}^{k}\times \mathbb{F}^{l}\times \mathbb{F}^{m}$};
\node (UMl) at (1.5,1.5) {$(\mathbb{F}^{k}\oplus\mathbb{F}^{l})\times\mathbb{F}^{m}$}; 
\node (UMr) at (4.5,1.5) {$\mathbb{F}^{k}\times(\mathbb{F}^{l}\oplus\mathbb{F}^{m})$};
\node (Ml) at (0,0) {$\mathbb{F}^{k+l}\times\mathbb{F}^{m}$};
\node (Mm) at (3,0) {$\mathbb{F}^{k}\oplus\mathbb{F}^{l}\oplus\mathbb{F}^{m}$};
\node (Mr) at (6,0) {$\mathbb{F}^{k}\times\mathbb{F}^{l+m}$};
\node (LMl) at (1.5,-1.5) {$\mathbb{F}^{k+l}\oplus\mathbb{F}^{m}$};
\node (LMr) at (4.5,-1.5) {$\mathbb{F}^{k}\oplus\mathbb{F}^{l+m}$};
\node (L) at (3,-3){$\mathbb{F}^{k+l+m}$};

\node (Label) at (3,-1.4){\tiny Lemma \ref{BVlemma}};

\path[->, font=\scriptsize,>=angle 90] 

(U) edge (UMl)
(U) edge (UMr)
(UMl) edge (Ml)
(UMl) edge (Mm)
(UMr) edge (Mm)
(UMr) edge (Mr)
(Ml) edge (LMl)
(Mm) edge (LMl)
(Mm) edge (LMr)
(Mr) edge (LMr)
(LMl) edge (L)
(LMr) edge (L);

\draw [->] (U) to [out=0, in=90] (Mr);
\draw [->] (U) to [out=180,in=90] (Ml);
\draw [->] (Mr) to [out=-90,in=0] (L);
\draw [->] (Ml) to [out=-90,in=180] (L);
\end{tikzpicture}
\end{equation}
\begin{equation}\label{CommDia2}
\begin{tikzpicture}[baseline=(current bounding box.center)]
\node (Lu) at (0,3){$\mathbb{F}^{k}\times\mathbb{F}^{l}$};
\node (Ru) at (4,3){$\mathbb{F}^{l}\times\mathbb{F}^{k}$};
\node (Lm) at (1,1.5){$\mathbb{F}^{k}\oplus\mathbb{F}^{l}$};
\node (Rm) at (3,1.5){$\mathbb{F}^{l}\oplus\mathbb{F}^{k}$};
\node (Ml) at (2,0){$\mathbb{F}^{k+l}$};
\node (Label) at (2,1.1){\tiny Lemma \ref{BVlemma}};

\path[->, font=\scriptsize,>=angle 90] 

(Lu) edge (Lm)
(Lu) edge node [above]{$\tau$}(Ru)
(Ru) edge (Rm)
(Lm) edge (Rm)
(Lm) edge (Ml)
(Rm) edge (Ml);

\draw [->] (Lu) to [out=-90,in=180] (Ml);
\draw [->] (Ru) to [out=-90,in=0] (Ml);

\end{tikzpicture}
\end{equation}
\end{center}

Similarly, one can define the $H$-space structure on $\mathbb{Z}\times G(\mathbb{F}^{\infty})$. We first fix an isometric isomorphism 
\[\alpha:\mathbb{F}^{\infty}\oplus\mathbb{F}^{\infty}\rightarrow\mathbb{F}^{\infty}\] 
and then define the (additive) operation to be the composition
\[ \{k\}\times G(\mathbb{F}^{\infty})\times \{l\}\times G(\mathbb{F}^{\infty})\rightarrow G(\mathbb{F}^{\infty}\oplus \mathbb{F}^{\infty})\rightarrow \{k+l\}\times G(\mathbb{F}^{\infty}),\]
where the first map is induced from the colimit of the commutative diagram below
\begin{center}
\begin{tikzpicture}
\node(Lu) at (0,2) {$G(\mathbb{F}^{k})\times G(\mathbb{F}^{l})$};
\node(Ll) at (0,0) {$G(\mathbb{F}^{k+1})\times G(\mathbb{F}^{l+1})$}; 
\node(Ru) at (5,2) {$G(\mathbb{F}^{k}\oplus\mathbb{F}^{l})$};
\node(Rl) at (5,0) {$G(\mathbb{F}^{k+1}\oplus \mathbb{F}^{l+1})$};

\path[->, font=\scriptsize,>=angle 90]

(Lu) edge (Ru)  
(Lu) edge node [right]{$(\iota_{k},\iota_{l})$}(Ll)
(Ll) edge (Rl) 
(Ru) edge (Rl);
\end{tikzpicture}
\end{center}
and, in the diagram above, the vertical map on the right-hand side is given as follows:
\begin{align*}
G(\mathbb{F}^{k}\oplus \mathbb{F}^{l})&\rightarrow G(\mathbb{F}^{k+1}\oplus \mathbb{F}^{l+1})\\
(V,W)&\mapsto V\oplus W\oplus e_{k+1}\otimes e_{1}\oplus e^{\prime}_{l+1}\otimes e_{1}
\end{align*}
where we have used the decomposition $(\mathbb{F}^{k+1}\oplus\mathbb{F}^{l+1})\otimes \mathbb{F}^{\infty}=\mathbb{F}^{k}\otimes\mathbb{F}^{\infty}\oplus\mathbb{F}^{l}\otimes\mathbb{F}^{\infty}\oplus\mathbb{F}\otimes\mathbb{F}^{\infty}\oplus\mathbb{F}\otimes\mathbb{F}^{\infty}$, and to distinguish, we let $<e_{i}^{\prime}>$ denote the standard basis on $\mathbb{F}^{l+1}$. The additive identity is $\ast=G(\mathbb{F}^{0})\hookrightarrow \{0\}\times G(\mathbb{F}^{\infty})$ by the homotopy commutative diagram
\begin{center}
\begin{tikzpicture}
\node(Lu) at (0,2) {$G(\mathbb{F}^{0}\oplus \mathbb{F}^{\infty})/G(\mathbb{F}^{\infty}\oplus \mathbb{F}^{0})$};
\node(Ru) at (7,2) {$G(\mathbb{F}^{\infty}\oplus \mathbb{F}^{\infty})$};
\node(Rl) at (7,0) {$G(\mathbb{F}^{\infty})$};
\path[->, font=\scriptsize,>=angle 90] 
(Lu) edge node [above]{$G(i\oplus \operatorname{id})/G(\operatorname{id}\oplus i)$}(Ru)  
(Lu) edge node [below]{$G(\pi)$}(Rl) 
(Ru) edge node [right]{$G(\alpha)$}(Rl);
\end{tikzpicture}
\end{center}
where $i$ is the inclusion $\mathbb{F}^{0}\hookrightarrow \mathbb{F}^{\infty}$ and $\pi:\mathbb{F}^{0}\oplus\mathbb{F}^{\infty} (\text{ resp. } \mathbb{F}^{\infty}\oplus\mathbb{F}^{0}) \rightarrow \mathbb{F}^{\infty}$ is the obvious isomorphism. Note that this homotopy in the diagram above might not preserve the base point, yet since $G(\mathbb{F}^{\infty})$ is a $\operatorname{CW}$-complex, one always can deform the operation such that it preserves the base point \cite[Remark $B.2.i$]{Wang5}.

Now, to see the associativity and commutativity of the additive operation, one only needs to replace $\mathbb{F}^{k}$ in the diagrams \eqref{AssoDia2} and \eqref{CommDia2} by $\mathbb{F}^{\infty}$. Lastly, the homotopy commutative diagram
\begin{center}
\begin{tikzpicture}
\node(Ll) at (0,0) {$G(\mathbb{F}^{\infty})\times G(\mathbb{F}^{\infty})$};
\node(Lu) at (0,2) {$G(\mathbb{F}^{k})\times G(\mathbb{F}^{l})$};
\node(Ml) at (4,0) {$G(\mathbb{F}^{\infty}\oplus\mathbb{F}^{\infty})$};
\node(Mu) at (4,2) {$G(\mathbb{F}^{k}\oplus\mathbb{F}^{l})$};
\node(Rl) at (8,0) {$G(\mathbb{F}^{\infty})$};
\node(Ru) at (8,2) {$G(\mathbb{F}^{k+l})$};

\node (Label) at (6,1){\tiny Lemma \ref{BVlemma}};

\path[->, font=\tiny,>=angle 90]

(Ll) edge (Ml)
(Ml) edge (Rl)
(Lu) edge (Mu)
(Mu) edge (Ru)  
(Lu) edge (Ll)
(Mu) edge (Ml)  
(Ru) edge (Rl);

\end{tikzpicture}
\end{center}
implies that the map
\[\coprod_{k}G(\mathbb{F}^{k})\rightarrow \mathbb{Z}\times G(\mathbb{F}^{\infty})\]
is a $H$-map, where the vertical maps are induced by the canonical inclusion $\mathbb{F}^{n}\hookrightarrow\mathbb{F}^{\infty}$. 
\end{proof}

Now we want to show these two $H$-spaces $\coprod\limits_{k}G(\mathbb{F}^{k})$ and $\coprod\limits_{k}\operatorname{Gr}_{k}(\mathbb{F}^{\infty})$ are homotopy equivalent as $H$-spaces. it is clear there is a canonical map 
\[\iota:\coprod_{k}\operatorname{Gr}_{k}(\mathbb{F}^{\infty})\rightarrow \coprod_{k}G(\mathbb{F}^{k})\]
given by the isometry 
\begin{align*}
\mathbb{F}^{\infty}=\mathbb{F}\otimes\mathbb{F}^{\infty}&\hookrightarrow\mathbb{F}^{k}\otimes\mathbb{F}^{\infty}\\
v=e_{1}\otimes v&\mapsto e_{1}\otimes v
\end{align*}
when $k>0$; for $k=0$, it is the identity since both of the spaces $\operatorname{Gr}_{0}(\mathbb{F}^{\infty})$ and $G(\mathbb{F}^{0})$ are a single point.

\begin{lemma}
The map 
\[\iota:\coprod_{k}\operatorname{Gr}_{k}(\mathbb{F}^{\infty})\rightarrow \coprod_{k}G(\mathbb{F}^{k})\]
is a homotopy equivalence of $H$-spaces.
\end{lemma}

\begin{proof}
It is clear that the additive identity of $\coprod\limits_{k}\operatorname{Gr}_{k}(\mathbb{F}^{\infty})$ is sent to the additive identity of $\coprod\limits_{k}G(\mathbb{F}^{k})$. To see the map $\iota$ respects $H$-structures, it suffices to note the homotopy commutative diagram below 
\begin{center}
\begin{tikzpicture}
\node(Ll) at (0,0) {$G(\mathbb{F}^{k})\times G(\mathbb{F}^{l})$};
\node(Lu) at (0,2) {$\operatorname{Gr}_{k}(\mathbb{F}^{\infty})\times \operatorname{Gr}_{l}(\mathbb{F}^{\infty})$};
\node(Ml) at (4,0) {$G(\mathbb{F}^{k}\oplus\mathbb{F}^{l})$};
\node(Mu) at (4,2) {$\operatorname{Gr}_{k+l}(\mathbb{F}^{\infty}\oplus\mathbb{F}^{\infty})$};
\node(Rl) at (8,0) {$G(\mathbb{F}^{k+l})$};
\node(Ru) at (8,2) {$\operatorname{Gr}_{k+l}(\mathbb{F}^{\infty})$};

\node (Label) at (6,1){\tiny Lemma \ref{BVlemma}};

\path[->, font=\tiny,>=angle 90]

(Ll) edge (Ml)
(Ml) edge (Rl)
(Lu) edge (Mu)
(Mu) edge (Ru)  
(Lu) edge (Ll)
(Mu) edge (Ml)  
(Ru) edge (Rl);

\end{tikzpicture}
\end{center}  
Finally, by lemma \ref{BVlemma}, any isomertic isomorphism $\mathbb{F}^{k}\otimes\mathbb{F}^{\infty}\rightarrow\mathbb{F}^{\infty}$ induces a homotopy inverse of $\iota$.
\end{proof}
 
Now, observe that the two sequences of spaces $\{G(\mathbb{F}^{k}),\iota_{k}\}$ and $\{\operatorname{Gr}_{k}(\mathbb{F}^{\infty}),\iota^{\prime}_{k}\}$ are compatible in the sense that the diagram 
\begin{center}
\begin{equation}\label{Stabilization}
\begin{tikzpicture}[baseline=(current bounding box.center)]
\node(Ll) at (0,0) {$G(\mathbb{F}^{k})$};
\node(Lu) at (0,2) {$\operatorname{Gr}_{k}(\mathbb{F}^{\infty})$};
\node(Ml) at (5,0) {$G(\mathbb{F}^{k}\oplus\mathbb{F})$};
\node(Mu) at (5,2) {$\operatorname{Gr}_{k+1}(\mathbb{F}^{\infty}\oplus\mathbb{F}^{\infty})$};
\node(Rl) at (8,0) {$G(\mathbb{F}^{k+1})$};
\node(Ru) at (8,2) {$\operatorname{Gr}_{k+1}(\mathbb{F}^{\infty})$};

\node (Label) at (6.5,1){\tiny Lemma \ref{BVlemma}};

\path[->, font=\tiny,>=angle 90]

(Lu) edge node [above]{$V\mapsto V\oplus e_{1}$}(Mu)
(Mu) edge (Ru)
(Ll) edge node [above]{$V\mapsto V\oplus e_{1}\otimes e_{1}$}(Ml)
(Ml) edge (Rl)  
(Lu) edge (Ll)
(Mu) edge (Ml)  
(Ru) edge (Rl);
\draw[->](Ll) to [out=-20,in=-160] node [above]{\scriptsize $\iota_{k}$}(Rl);
\draw[->](Lu) to [out=20,in=160] node [below]{\scriptsize $\iota^{\prime}_{k}$} (Ru);
\end{tikzpicture}
\end{equation}
\end{center}
commutes up to homotopy. Applying the (ordinary) homology functor $H_{\ast}$, the diagram \eqref{Stabilization} induces the following commutative diagram of graded abelian groups  
\begin{center}
\begin{tikzpicture}
\node(Ll) at (0,0) {$H_{\ast}(G(\mathbb{F}^{k}))$};
\node(Lu) at (0,2) {$H_{\ast}(\operatorname{Gr}_{k}(\mathbb{F}^{\infty}))$};
\node(Ml) at (4,0) {$H_{\ast}(G(\mathbb{F}^{k+1}))$};
\node(Mu) at (4,2) {$H_{\ast}(\operatorname{Gr}_{k+1}(\mathbb{F}^{\infty}))$};
\node(Rl) at (8,0) {$H_{\ast}(G(\mathbb{F}^{k+2}))$};
\node(Ru) at (8,2) {$H_{\ast}(\operatorname{Gr}_{k+2}(\mathbb{F}^{\infty}))$};
 
\path[->, font=\tiny,>=angle 90]

(Lu) edge node [above]{$\iota_{k,\ast}^{\prime}$}(Mu)
(Mu) edge node [above]{$\iota_{k+1,\ast}^{\prime}$}(Ru)
(Ll) edge node [above]{$\iota_{k,\ast}$}(Ml)
(Ml) edge node [above]{$\iota_{k+1,\ast}$}(Rl)  
(Lu) edge (Ll)
(Mu) edge (Ml)  
(Ru) edge (Rl);

\end{tikzpicture}
\end{center}
Since the stabilization maps are induced by adding an element in $\operatorname{Gr}_{1}(\mathbb{F}^{\infty})$ (resp. $ G(\mathbb{F}^{1})$) that generates $\pi_{0}(\coprod_{k}\operatorname{Gr}_{k}(\mathbb{F}^{\infty}))$ (resp. $\pi_{0}(\coprod_{k}G(\mathbb{F}^{k}))$), we obtain an isomorphism of rings
\begin{equation}\label{Isombetweentwomodels}
H_{\ast}(\coprod_{k}\operatorname{Gr}_{k}(\mathbb{F}^{\infty}),R)[\pi_{0}^{-1}]\xrightarrow{\sim} H_{\ast}(\coprod
_{k}G(\mathbb{F}^{k}),R)[\pi_{0}^{-1}],
\end{equation} 
where $\pi_{0}$ stands for the abelian monoid of connect components in $\coprod\limits_{k}\operatorname{Gr}_{k}(\mathbb{F}^{\infty})$ (resp. $\coprod\limits_{k}G(\mathbb{F}^{k})$) and $R$ any commutative ring.

\begin{remark}
In the argument above, choosing which element in $\operatorname{Gr}_{1}(\mathbb{F}^{\infty})$ (resp. $G_{1}(\mathbb{F})$) is not really important because any two of them give homotopic maps from $\operatorname{Gr}_{k}(\mathbb{F}^{\infty})$ to $\operatorname{Gr}_{k+1}(\mathbb{F}^{\infty})$ (resp. $G(\mathbb{F}^{k})$ to $G(\mathbb{F}^{k+1})$).
\end{remark}

We are now in the position to state the main theorem of this section.

\begin{theorem}\label{GroupcomplofG}
The composition 
\[\coprod_{k}\operatorname{Gr}_{k}(\mathbb{F}^{\infty})\rightarrow \coprod_{k}G(\mathbb{F}^{k})\rightarrow \mathbb{Z}\times G(\mathbb{F}^{\infty})\]
is a topological group completion. 
\end{theorem}
\begin{proof}
Recall that, given two $H$-spaces $X$, $Y$, a map $f:X\rightarrow Y$ is a topological group completion if and only if the following three statements are satisfied: firstly, $f$ is a $H$-map; secondly, $Y$ is group-like; thirdly, the induced homomorphism $H_{\ast}(X,R)[\pi_{0}(X)^{-1}]\xrightarrow{f_{\ast}} H_{\ast}(Y,R)$ is an isomorphism of rings, for any commutative ring $R$.

The first condition is easy as both maps in the statement are $H$-maps. It is also clear that $\pi_{0}(\mathbb{Z}\times G(\mathbb{F}^{\infty}))=\mathbb{Z}$, and hence, the second condition is satisfied as well. As for the third condition, in view of the isomorphism \eqref{Isombetweentwomodels}, it suffices to show the homomorphism   
\[H_{\ast}(\coprod_{k}G(\mathbb{F}^{k}),R)[\pi_{0}^{-1}]\xrightarrow{\sim} H_{\ast}(\mathbb{Z}\times G(\mathbb{F}^{\infty}),R)\]
is an isomorphism of rings.
To see this, we note first that, for any sequence of cofibrations, the colimit commutes with the homology functor, and hence the homomorphism of $H_{\ast}(\coprod\limits_{k}G(\mathbb{F}^{k}),R)$-modules   
\begin{equation}\label{IsoofModules}
H_{\ast}(\coprod_{k}G(\mathbb{F}^{k}),R)[\pi_{0}^{-1}]\xrightarrow{\sim} H_{\ast}(\mathbb{Z}\times
G(\mathbb{F}^{\infty}),R)
\end{equation}
is in fact an isomorphism. Then, by the homotopy commutative diagram  
\begin{center} 
\begin{tikzpicture}
\node(Ll) at (0,0) {$\{k-m\}\times G(\mathbb{F}^{\infty})\times \{l-n\}\times G(\mathbb{F}^{\infty})$};
\node(Lu) at (0,1.5) {$(G(\mathbb{F}^{k}),m)\times (G(\mathbb{F}^{l}),n)$};
\node(Ml) at (4.5,-1) {$G(\mathbb{F}^{\infty}\oplus\mathbb{F}^{\infty})$};
\node(Mu) at (4.5,2.5) {$G(\mathbb{F}^{k}\oplus\mathbb{F}^{l})$};
\node(Rl) at (8,0) {$\{(k+l)-(m+n)\}\times G(\mathbb{F}^{\infty})$};
\node(Ru) at (8,1.5) {$(G(\mathbb{F}^{k+l}),m+n)$};

\node (Label) at (6,.75){\tiny Lemma \ref{BVlemma}};

\path[->, font=\tiny,>=angle 90]

(Lu) edge (Ll)
(Mu) edge (Ml)  
(Ru) edge (Rl);
\draw[->] (Ll) to [out=-40,in=180] (Ml);
\draw[->] (Ml) to [out=0,in=-140] (Rl);
\draw[->] (Lu) to [out=40,in=-180] (Mu);
\draw[->] (Mu) to [out=0,in=140] (Ru);
\end{tikzpicture}
\end{center} 
we see the multiplicative structures are also preserved, where $m,n$ indicate that they are the $m$-th and $n$-th copies of $G(\mathbb{F}^{k})$ and $G(\mathbb{F}^{l})$ in 
\[\{ \coprod_{k}G(\mathbb{F}^{k}),\iota:=\coprod_{k}\iota_{k} \},\]
respectively.
\end{proof}

It is well-known that, for every compact Hausdorff space $X$, there are isomorphisms of abelian monoids and abelian groups: 
\begin{align*}
\coprod_{k}\operatorname{Vect}^{\mathbb{F}}_{k}(X)\cong [X,\coprod_{k}\operatorname{Gr}_{k}(\mathbb{F}^{\infty})];\\
K^{\mathbb{F}}(X)\cong [X,\mathbb{Z}\times G(\mathbb{F}^{\infty})],
\end{align*}
where $\operatorname{Vect}^{\mathbb{R}}_{k}(X)$ (resp. $\operatorname{Vect}^{\mathbb{C}}_{k}(X)$) is the set of isomorphism classes of real (resp. complex) vector bundles over $X$ and $K^{\mathbb{R}}(X)$ (resp. $K^{\mathbb{C}}(X)$) is the real (resp. complex) topological $K$-theory of $X$. Combining this fact with Theorem \ref{GroupCompletionthesection}, we see the topological group completion
\[\coprod_{k}\operatorname{Gr}_{k}(\mathbb{F}^{\infty})\rightarrow \mathbb{Z}\times G(\mathbb{F}^{\infty})\]
realizes the algebraic group completion
\[\coprod_{k}\operatorname{Vect}^{\mathbb{F}}_{k}(X)\rightarrow K^{\mathbb{F}}(X).\] 
In particular, the Bott isomorphism of reduced $K$-groups   
\[\tilde{K}^{\mathbb{F}}(X)\xrightarrow{\sim} \tilde{K}^{\mathbb{F}}(X\wedge S^{d})\]
induces a map
\begin{equation}\label{Bottmapspacelevel}
B:\mathbb{Z}\times G(\mathbb{F}^{\infty})\rightarrow \Omega^{d}(\mathbb{Z}\times G(\mathbb{F}^{\infty})),
\end{equation} 
where $X$ is any pointed compact Hausdorff space and $d=2$ (resp. $d$=8) if $\mathbb{F}=\mathbb{C}$ (resp. $\mathbb{F}=\mathbb{R}$).

A priori, the map \eqref{Bottmapspacelevel} is determined only up to weak homotopy. However, in the case of complex topological $K$-theory, the Atiyah-Hirzebruch spectral sequence tells us the $\operatorname{lim}^{1}$-term vanishes \cite[p.301]{Sw}, namely
\[\operatorname*{lim^{1}}_{n}\tilde{K}^{\mathbb{C}}(\operatorname{Gr}_{n}(\mathbb{C}^{n}\otimes \mathbb{C}^{n})\wedge S^{2})=0=\operatorname*{lim^{1}}_{n}\tilde{K}^{\mathbb{C}}(\operatorname{Gr}_{n}(\mathbb{C}^{n}\otimes \mathbb{C}^{n})_{+}\wedge S^{2}),\]
and hence, the following homomorphisms 
\begin{align*}
\tilde{K}^{\mathbb{C}}(G(\mathbb{C}^{\infty})_{+}\wedge S^{2})&\xrightarrow{\sim} \operatorname*{lim}_{n}\tilde{K}^{\mathbb{C}}(\operatorname{Gr}_{n}(\mathbb{C}^{n}\otimes\mathbb{C}^{n})_{+}\wedge S^{2});\\
\tilde{K}^{\mathbb{C}}(G(\mathbb{C}^{\infty}) \wedge S^{2})&\xrightarrow{\sim} \operatorname*{lim}_{n}\tilde{K}^{\mathbb{C}}(\operatorname{Gr}_{n}(\mathbb{C}^{n}\otimes\mathbb{C}^{n}) \wedge S^{2})
\end{align*}
are isomorphisms. It implies the map \eqref{Bottmapspacelevel} is uniquely determined, up to homotopy, in the case of the complex numbers.

\begin{remark}
In general, given any compact Hausdorff space $A$, a topological group completion $X\rightarrow Y$ does not induce an algebraic group completion $[A,X]\rightarrow [A,Y]$, for example, algebraic $K$-theory of a ring. 
\end{remark}

\section{The standard category of vector spaces}
\subsection{The category $\mathcal{V}_{\mathbb{F}}$}
Recall that   the standard category of vector spaces over $\mathbb{F}$, denoted by $\mathcal{V}_{\mathbb{F}}$ is given by   
\begin{align*}
O(\mathcal{V}_{\mathbb{F}}):=\{\mathbb{F}^{i}\}_{i\in\mathbb{N}\cup\{0\}}
M(\mathcal{V}_{\mathbb{F}}):=\coprod_{m,n}M_{m,n}(\mathbb{F})
\end{align*}
where $M_{m,n}(\mathbb{F})$ is the space of $m$-by-$n$ matrices topologized as $\mathbb{F}^{mn}$.

The category $\mathcal{V}_{\mathbb{F}}$ is an enriched category, and it admits two functors. The first one, given by the direct sum of vector spaces, is defined by the assignment 
\begin{align}\label{DirectsumVf}   
\oplus:\mathcal{V}_{\mathbb{F}}\times \mathcal{V}_{\mathbb{F}}&\rightarrow \mathcal{V}_{\mathbb{F}}\nonumber\\
(\mathbb{F}^{n},\mathbb{F}^{m})&\mapsto \mathbb{F}^{n+m}\\
 (\{e_{i}\}_{i=1}^{n},\{e^{\prime}_{j}\}_{j=1}^{m})&\mapsto \{f_{k}=\begin{cases} (e_{k},0) & k\leq n\\
              (0,e^{\prime}_{k-n}) & i>n\end{cases} \hspace{2em} \},\nonumber
\end{align}
while, the second one, induced by the tensor produce of vector spaces, is given by the assignment  
\begin{align}\label{TensorVf} 
\otimes:\mathcal{V}_{\mathbb{F}}\times \mathcal{V}_{\mathbb{F}}&\rightarrow \mathcal{V}_{\mathbb{F}}\nonumber\\
(\mathbb{F}^{n},\mathbb{F}^{m})&\mapsto \mathbb{F}^{nm}\\
 (\{e_{i}\}_{i=1}^{n},\{e^{\prime}_{j}\}_{j=1}^{m})&\mapsto \{f_{(m-1)i+j}=e_{i}\otimes e^{\prime}_{j}\}\nonumber
\end{align}
These two functors induce two operations on $\vert iso\mathcal{V}_{\mathbb{F}}\vert$, still denoted by $\oplus$ and $\otimes$, 
\begin{align*} 
\oplus:\vert iso\mathcal{V}_{\mathbb{F}}\vert\times \vert iso\mathcal{V}_{\mathbb{F}}\vert&\rightarrow \vert iso\mathcal{V}_{\mathbb{F}}\vert;\\
\otimes: \vert iso\mathcal{V}_{\mathbb{F}}\vert\times \vert iso\mathcal{V}_{\mathbb{F}}\vert&\rightarrow \vert iso\mathcal{V}_{\mathbb{F}}\vert
\end{align*}
which give a $H$-semiring structure on the space $\vert iso\mathcal{V}_{\mathbb{F}}\vert$. In fact, via May's delooping machine, these two operations give the space $\vert iso\mathcal{V}_{\mathbb{F}}\vert$ an $E_{\infty}$-ring space structure \cite[Chapter VI and VII]{Ma4} (see also \cite{Ma7} and \cite{Ma8}), we are not pursuing higher structures in the present paper, however.  

\begin{lemma}
The operations $\oplus$, $\otimes$ induce a $H$-semiring structure on $\vert iso\mathcal{V}_{\mathbb{F}}\vert$ with $\mathbb{F}^{0}=0$ the additive identity and $\mathbb{F}$ the multiplicative identity.
\end{lemma}
\begin{proof}
Firstly, we note that, from the assignment \eqref{DirectsumVf} (resp. \eqref{TensorVf}), it is clear that $0$ (resp. $\mathbb{F}$) is the strict additive (resp. multiplicative) identity and the operation $\oplus$ (resp. $\otimes$) is associative. 

Now, to see the commutativity of the operation $\oplus$ (resp. $\otimes$), we observe that the diagram  
\begin{center}
\begin{tikzpicture}
\node(Lu) at (0,2) {$iso\mathcal{V}_{\mathbb{F}} \times  iso\mathcal{V}_{\mathbb{F}}$};
\node(Ru) at (6,2) {$iso\mathcal{V}_{\mathbb{F}} \times  iso\mathcal{V}_{\mathbb{F}}$};
\node(Ml) at (3,0) {$ iso\mathcal{V}_{\mathbb{F}} $};
\path[->, font=\scriptsize,>=angle 90] 
(Lu) edge node [above]{$\tau$} (Ru)  
(Lu) edge node [left, yshift=-.5em]{$\otimes/\oplus$}(Ml) 
(Ru) edge node [right, yshift=-.5em]{$\otimes/\oplus$}(Ml);
\end{tikzpicture}
\end{center}
where $\tau(\mathbb{F}^{n},\mathbb{F}^{m}):=(\mathbb{F}^{m},\mathbb{F}^{n})$, though not commutative on the nose, does commute up to natural transformations. For instance, in the case of the operation $\oplus$, the natural transformation is given by
\[(\mathbb{F}^{n},\mathbb{F}^{m})\mapsto 
\begin{bmatrix}
0&\operatorname{id}_{n}\\
\operatorname{id}_{m}&0
\end{bmatrix}.\]

Similarly, one can check the operation $\otimes$ is right- and left-distributive. The right distributive law follows from the commutative diagram 
\begin{center}
\begin{tikzpicture}
\node(Lu) at (0,2) {$iso\mathcal{V}_{\mathbb{F}}\times iso\mathcal{V}_{\mathbb{F}}\times iso\mathcal{V}_{\mathbb{F}}$};
\node(Lm) at (0,1) {$iso\mathcal
V_{\mathbb{F}}\times iso\mathcal{V}_{\mathbb{F}}\times iso\mathcal{V}_{\mathbb{F}}\times iso\mathcal{V}_{\mathbb{F}}$};
\node(Ll) at (0,0) {$iso\mathcal{V}_{\mathbb{F}}\times iso\mathcal{V}_{\mathbb{F}}$}; 
\node(Ru) at (5,2) {$iso\mathcal{V}_{\mathbb{F}}\times iso\mathcal{V}_{\mathbb{F}}$};
 
\node(Rl) at (5,0) {$iso\mathcal{V}_{\mathbb{F}}$};

\path[->, font=\scriptsize,>=angle 90]

(Lu) edge node [above]{$(\oplus,\operatorname{id})$}(Ru)  
(Lu) edge node [right]{$r$}(Lm)
(Lm) edge node [right]{$(\otimes,\otimes)$}(Ll)
(Ll) edge node [above]{$\oplus$}(Rl) 
(Ru) edge node [right]{$\otimes$}(Rl);

\end{tikzpicture}
\end{center} 
where $r(\mathbb{F}^{m},\mathbb{F}^{n},\mathbb{F}^{k}):=(\mathbb{F}^{m},\mathbb{F}^{k},\mathbb{F}^{n},\mathbb{F}^{k})$; The left distributive law follows from the same diagram with $r$ replaced by $l(\mathbb{F}^{m},\mathbb{F}^{n},\mathbb{F}^{k}):=(\mathbb{F}^{m},\mathbb{F}^{n},\mathbb{F}^{m},\mathbb{F}^{k})$ and $(\oplus, \operatorname{id})$ by $(\operatorname{id},\oplus)$. However, in this case, it commutes only up to natural transformations. A precise natural transformation is given by the assignment
\[(\mathbb{F}^{m},\mathbb{F}^{n},\mathbb{F}^{k})\mapsto
\begin{bmatrix}
\operatorname{id}_{n}&0&0&0&&..&0\\
  0&0&\operatorname{id}_{n}&0&&...&0\\
  0&0&0&0&\operatorname{id}_{n}&...&0\\
   &&\vdots&&&&\\
   0&\operatorname{id}_{k}&0&0&0&...&0\\
   0&0&0&\operatorname{id}_{k}&0&...&0\\
   &&\vdots&&&&\\
   0&0&0&0&0&...&\operatorname{id}_{k} 
\end{bmatrix}
\] 
\end{proof}
\subsection{Geometric meaning of $\mathcal{V}_{\mathbb{F}}$}
In view of the classification theorem \cite[Theorem $3.2$]{Wang4}, there is a $1$-$1$ correspondence
\begin{equation}\label{Vfcocyclesandhtyclasses}
[X,B iso\mathcal{V}_{\mathbb{F}}]\leftrightarrow iso\mathcal{V}_{\mathbb{F}}(X),
\end{equation}
for any paracompact Hausdorff space $X$, where $iso\mathcal{V}_{\mathbb{F}}(X)$ is the set of isomorphism classes of $iso\mathcal{V}_{\mathbb{F}}$-cocycles on $X$ (see \cite[Section I]{Ha} and \cite[Section $2$]{Wang4}). On the other hand, it is well-known that the set of isomorphism classes of $iso\mathcal{V}_{\mathbb{F}}$-cocycles on $X$ is the same as the set of isomorphism classes of vector bundles over $X$, denoted by $\coprod\limits_{k}\operatorname{Vect}_{k}(X).$ Now, since $\operatorname{id}_{n}\in\operatorname{GL}_{n}(\mathbb{F})$ is well-pointed, the simplicial space $\operatorname{Ner}_{\cdot}iso\mathcal{V}_{\mathbb{F}}$ is proper; combing this with Theorem $4.3$ in \cite{Wang4} and the fact that $\operatorname{Ner}_{i}iso\mathcal{V}_{\mathbb{F}}$ has the homotopy type of a $\operatorname{CW}$-complex, for every $i$, we see the map 
\[Bsio\mathcal{V}_{\mathbb{F}}\rightarrow \vert\operatorname{Ner}_{\cdot}iso\mathcal{V}_{\mathbb{F}}\vert\]
is a homotopy equivalence, and thus we recover the well-known result.  
\begin{theorem}\label{GeomeaningofVf}
There is a $1$-$1$ correspondence
\[[X,\vert iso\mathcal{V}_{\mathbb{F}}\vert]\leftrightarrow \coprod_{k}\operatorname{Vect}^{\mathbb{F}}_{k}(X).\] 
\end{theorem}
\section{A new category of vector spaces} 
\subsection{The category $\mathcal{V}_{\mathbb{F}}^{f}$}
Before constructing the fat category of vector spaces over $\mathbb{F}$, denoted by $\mathcal{V}_{\mathbb{F}}^{f}$, we recall the manifold structure of $\operatorname{Gr}_{k}(\mathbb{F}^{m})$, the Grassmannian of $k$-dimensional subspaces in $\mathbb{F}^{m}$.   

The standard coordinate charts of the Grassmannian
$\operatorname{Gr}_{k}(\mathbb{F}^{m})$ can be described as follows: Given a point $V\in \operatorname{Gr}_{k}(\mathbb{F}^{m})$, one has the associated projection $\pi_{V}:\mathbb{F}^{m}\rightarrow V$ defined by 
\begin{align*}
\pi_{V}\vert_{V}&=\operatorname{id};\\
\pi_{V}\vert_{V^{\perp}}&=0,
\end{align*}  
where $V^{\perp}$ is the orthogonal complement of $V$ with respect to the standard inner product of $\mathbb{F}^{m}$. Using this projection $\pi_{V}$, one can define a bijection $\phi_{V}$ from $\mathfrak{N}_{V}$, the set of all $k$-planes $W$ with $\pi_{V}(W)=V$, to $M_{m-k,k}(\mathbb{F})$, the space of $(m-k)$-by-$(k)$ matrices, which is topologized as $\mathbb{F}^{(m-k)k}$. In more details, we first choose an orthogonal basis of $\mathbb{F}^{m}$ that contains a basis of $V$, denoted by $v$, and a basis of $V^{\perp}$, denoted by $v^{\perp}$, then the matrix $A$ associated to $W\in \mathbb{F}^{m}$ is given by the identity 
\[<\mathbf{v}+\mathbf{v}^{\perp}A>=W.\] 
Furthermore, it is not difficult check, given two elements $V,V^{\prime}\in \operatorname{Gr}_{k}(\mathbb{F}^{m})$, the map  
\[\phi_{V}\circ\phi^{-1}_{V^{\prime}}\vert_{\phi_{V^{\prime}}(\mathfrak{N}_{V}\cap\mathfrak{N}_{V^{\prime}})}:\phi_{V^{\prime}}(\mathfrak{N}_{V}\cap\mathfrak{N}_{V^{\prime}})\rightarrow M_{m-k,k}(\mathbb{F})\]
is a diffeomorphism whenever it can be defined. More precisely, given $B\in \phi_{V^{\prime}}(\mathfrak{N}_{V}\cap\mathfrak{N}_{V^{\prime}})$, the matrix $\phi_{V}\circ\phi^{-1}_{V^{\prime}}(B)$ can be written in terms of $B$ and the transition matrix between the bases $<v,v^{\perp}>$ and $<v^{\prime},v^{\prime,\perp}>$. Hence, the collection $\{\mathfrak{N}_{V},\phi_{V}\}$ comprises a smooth structure on $\operatorname{Gr}_{k}(\mathbb{F}^{m})$. In particular, the canonical inclusion $\mathbb{F}^{m}\xrightarrow{x\mapsto (x,0)} \mathbb{F}^{m+1}$ gives an smooth embedding
\begin{equation}\label{cofibbetGr}
\operatorname{Gr}_{k}(\mathbb{F}^{m})\hookrightarrow \operatorname{Gr}_{k}(\mathbb{F}^{m+1})
\end{equation}  
and hence a cofibration.  

Now, observe that there is a vector bundle 
\[(s^{m}_{k},t^{n}_{l}):\operatorname{Mor}_{k,l}^{m,n}(\mathbb{F})\rightarrow \operatorname{Gr}_{k}(\mathbb{F}^{m})\times \operatorname{Gr}_{l}(\mathbb{F}^{n})\]
whose total space is the set of linear maps from $V$ to $W$, for every $V\in \operatorname{Gr}_{k}(\mathbb{F}^{m})$ and $W\in \operatorname{Gr}_{l}(\mathbb{F}^{n})$. We topologize it in the following way: Given $X_{0}\in \operatorname{Gr}_{k}(\mathbb{F}^{m})$ and  $Y_{0}\in \operatorname{Gr}_{l}(\mathbb{F}^{n})$, we define $\mathfrak{N}_{X_{0},Y_{0}}:= (s^{m}_{k},t^{n}_{l})^{-1}(\mathfrak{N}_{Y_{0}}\times\mathfrak{N}_{Y_{0}})$ and observe that there is a bijection
\begin{align} 
\phi_{X_{0},Y_{0}}:\mathfrak{N}_{X_{0},Y_{0}}&\rightarrow 
M_{m-k,k}(\mathbb{F})\times M_{n-l,l}(\mathbb{F})\times M_{l,k}(\mathbb{F}).\label{CoordinateMor}\\
(T:X\rightarrow Y)&\mapsto (A_{X},B_{Y},[T])\nonumber
\end{align}
where $A_{X}$, $B_{Y}$ are the matrices that satisfy 
\begin{align*}
<\mathbf{x}=\mathbf{x}_{0}+\mathbf{x}_{0}^{\perp}A_{X}>&=X,\\
<\mathbf{y}=\mathbf{y}_{0}+\mathbf{y}_{0}^{\perp}A_{Y}>&=Y,
\end{align*}
respectively, and $[T]$ is the matrix that realizes $T$ in terms of the bases $\mathbf{x}$ and $\mathbf{y}$. In addition, given $X_{0},X_{1}\in \operatorname{Gr}_{k}(\mathbb{F}^{m})$ and $Y_{0},Y_{1}\in \operatorname{Gr}_{l}(\mathbb{F}^{n})$, whenever the map 
\[\phi_{X_{0},Y_{0}}\circ\phi_{X_{1},Y_{1}}^{-1}\vert_{\phi_{X_{0},Y_{0}}(\mathfrak{N}_{X_{0},Y_{0}}\cap\mathfrak{N}_{X_{1},Y_{1}})}\] 
is defined, the $l$-by-$k$ matrix in $\phi_{X_{0},Y_{0}}\circ\phi_{X_{1},Y_{1}}^{-1}(A_{X},B_{Y},[T])$ can be expressed as $C[S]D$, where $D$ and $C$ are the $k$-by-$k$ matrix and the $l$-by-$l$ matrix that are induced from the transition matrix between $<x_{0},x_{0}^{\perp}>$ and $<x_{1},x_{1}^{\perp}>$ and the one between $<y_{0},y_{0}^{\perp}>$ and $<y_{1},y_{1}^{\perp}>$. Thus, we obtain a chart of the vector bundle 
\[(s^{m}_{k},t^{n}_{l}):\operatorname{Mor}_{k,l}^{m,n}(\mathbb{F})\rightarrow \operatorname{Gr}_{k}(\mathbb{F}^{m})\times \operatorname{Gr}_{l}(\mathbb{F}^{n}),\]
and hence it is topologized. 

Now, observe that the canonical inclusions $\mathbb{F}^{m}\hookrightarrow \mathbb{F}^{m+1}$ and $\mathbb{F}^{n}\hookrightarrow \mathbb{F}^{n+1}$ give us the following two maps
\begin{align}
\operatorname{Mor}_{k,l}^{m,n}(\mathbb{F})&\rightarrow \operatorname{Mor}_{k,l}^{m+1,n}(\mathbb{F})\nonumber;\\
\operatorname{Mor}_{k,l}^{m,n}(\mathbb{F})&\rightarrow \operatorname{Mor}_{k,l}^{m,n+1}(\mathbb{F})\label{cofibbetMor},
\end{align}
which are, in fact, cofibrations, in view of the pullback diagram below and Lemma \ref{PullbackandCof}
\begin{center}
\begin{tikzpicture}
\node(Lu) at (0,2) {$\operatorname{Mor}^{m,n}_{k,k}(\mathbb{F})$};
\node(Ll) at (0,0) {$\operatorname{Gr}_{k}(\mathbb{F}^{m})\times \operatorname{Gr}_{k}(\mathbb{F}^{n})$}; 
\node(Ru) at (5,2) {$\operatorname{Mor}^{m+1,n+1}_{k,k}(\mathbb{F})$};
\node(Rl) at (5,0) {$\operatorname{Gr}_{k}(\mathbb{F}^{m+1})\times \operatorname{Gr}_{k}(\mathbb{F}^{n+1})$};

\path[->, font=\scriptsize,>=angle 90]

(Lu) edge (Ru)  
(Lu) edge (Ll)
(Ll) edge (Rl) 
(Ru) edge (Rl);
\end{tikzpicture}
\end{center}  

With these preliminaries in place, we can now define the category $\mathcal{V}^{f}_{\mathbb{F}}$. The space of objects and the space of morphisms of this internal category in $\mathpzc{Top}^{w}$ are given by the colimits 
\begin{align}
O(\mathcal{V}_{\mathbb{F}}^{f})&:= \coprod_{k}\operatorname*{colim}_{m} \operatorname{Gr}_{k}(\mathbb{F}^{m})\nonumber\\
M(\mathcal{V}_{\mathbb{F}}^{f})&:= \coprod_{k,l}\operatorname*{colim}_{m,n} \operatorname{Mor}_{k,l}^{m,n}(\mathbb{F})\label{colimitofMor},
\end{align}
respectively, where the colimits are taken over the systems of cofibrations given in \eqref{cofibbetGr} and \eqref{cofibbetMor}, respectively.

\begin{lemma}\label{InternalCatVft}
$\mathcal{V}_{\mathbb{F}}^{f}$ is an internal category in $\mathpzc{Top}^{w}$.
\end{lemma} 
\begin{proof}
The source and target maps, denoted by $s$ and $t$, are given by
the colimits
\begin{align*} 
s:= \coprod_{k}\operatorname*{colim}_{m}s_{k}^{m},\\ 
t:= \coprod_{l}\operatorname*{colim}_{m}t_{k}^{m};
\end{align*}
and the identity-assigning map  
\[e:= \coprod_{k}\operatorname*{colim}_{m}e_{k}^{m},\]
where $e^{m}_{k}$ is the map
\begin{align*}
e_{k}^{m}:\operatorname{Gr}_{k}(\mathbb{F}^{m})&\rightarrow \operatorname{Mor}_{k,k}^{m,m}(\mathbb{F})\\ 
V&\mapsto \operatorname{id}:V\rightarrow V.  
\end{align*} 
Now the composition of matrices induces a map
\begin{equation}\label{FinitecompositionforVft}
\operatorname{Mor}_{k,l}^{m,n}(\mathbb{F})\times_{\operatorname{Gr}_{l}(\mathbb{F}^{n})}\operatorname{Mor}_{l,j}^{n,p}(\mathbb{F})
\rightarrow \operatorname{Mor}_{k,j}^{m,p}(\mathbb{F}),
\end{equation}
for every $m,n,p,k,l,j\in \mathbb{N}\cup \{0\}$. To see its continuity, we check its restriction to the coordinate charts \eqref{CoordinateMor}---for simplicity, $\mathbb{F}$ is dropped from the notation:
\begin{align*}
\mathfrak{N}_{T_{0}:X_{0}\rightarrow Y_{0}}&\cong M_{m-k,k}\times M_{n-l,l}\times M_{l,k},\\ 
\mathfrak{N}_{S_{0}:Y_{0}\rightarrow Z_{0}}&\cong M_{n-l,l}\times M_{p-j,j}\times M_{j,l}.
\end{align*}
Then the composition \eqref{FinitecompositionforVft} can be realized in terms of these charts as follows:
\begin{align*}
(M_{m-k,k} \times M_{n-l,l} \times M_{l,k}) \times_{M_{n-l,l}}(M_{n-l,l}\times M_{p-j,j}\times M_{j,l}) 
&\\
\longrightarrow  M_{m-k,k}&\times M_{p-j,j}\times M_{j,k},\\
(A_{X},B_{Y},[T];B_{Y},C_{Z},[S])\mapsto (A_{X},C_{Z}&,[S]\circ[T]).
\end{align*} 
Hence, the composition \eqref{FinitecompositionforVft} is continuous, and it gives us a map 
\[\operatorname*{colim}_{m}(\operatorname{Mor}_{k,l}^{m,m}(\mathbb{F})\times_{\operatorname{Gr}_{l}(\mathbb{F}^{m})}\operatorname{Mor}_{l,j}^{m,m}(\mathbb{F}))\rightarrow \operatorname*{colim}_{m}\operatorname{Mor}_{k,j}^{m,m}(\mathbb{F}).\]
Employing Lemma \ref{Colimitandpullback}, which implies the homeomorphism  
\begin{multline*}
\operatorname*{colim}_{m}(\operatorname{Mor}_{k,l}^{m,m}(\mathbb{F})\times_{\operatorname{Gr}_{l}(\mathbb{F}^{m})}\operatorname{Mor}_{l,j}^{m,m}(\mathbb{F}))\\
=\operatorname*{colim}_{m}\operatorname{Mor}_{k,l}^{m,m}(\mathbb{F})\times_{\operatorname*{colim}_{m}\operatorname{Gr}_{l}(\mathbb{F}^{m})}\operatorname*{colim}\limits_{m}\operatorname{Mor}_{l,j}^{m,m}(\mathbb{F}),
\end{multline*}
we obtain the required composition 
\[\circ:M(\mathcal{V}_{\mathbb{F}}^{f})\times_{O(\mathcal{V}_{\mathbb{F}}^{f})}M(\mathcal{V}_{\mathbb{F}}^{f})\rightarrow M(\mathcal{V}_{\mathbb{F}}^{f})\]
in $\mathcal{V}_{\mathbb{F}}^{f}$. The compatibility of these structure maps can be derived directly from the definitions. 
\end{proof}
 
Like $\mathcal{V}_{\mathbb{F}}$, the category $\mathcal{V}_{\mathbb{F}}^{f}$ admits two operations given by the direct sum and tensor product of vector spaces. The construction is similar to those in $\mathcal{V}_{\mathbb{F}}$ but more delicate.  

Recall first, given a vector space $V$, we have the manifold of $k$-dimensional subspace of $V$, denoted by $\operatorname{Gr}_{k}(V)$, and an isomorphism of vector spaces $\phi :V\xrightarrow{\sim} W$ induces a diffeomorphism  
\[\phi_{\ast}:\operatorname{Gr}_{k}(V)\rightarrow \operatorname{Gr}_{k}(W).\]
Now consider the composition
\begin{equation}\label{OplusforOVFf}
\operatorname{Gr}_{k}(\mathbb{F}^{n})\times\operatorname{Gr}_{k^{\prime}}(\mathbb{F}^{n})\rightarrow\operatorname{Gr}_{k+k^{\prime}}(\mathbb{F}^{n}\oplus\mathbb{F}^{n})\rightarrow \operatorname{Gr}_{k+k^{\prime}}(\mathbb{F}^{2n}),
\end{equation}
where the first map sends the pair $(X,Y)$ to $X\oplus Y\in \operatorname{Gr}_{k+l}(\mathbb{F}^{n}\oplus \mathbb{F}^{n})$ and the second map sends $X\oplus Y$ to $\theta_{\ast}(X\oplus Y)\in \operatorname{Gr}_{k+k^{\prime}}(\mathbb{F}^{2n})$ via the isomorphism 
\begin{align} 
\theta:\mathbb{F}^{n}\oplus \mathbb{F}^{n}&\xrightarrow{\sim} \mathbb{F}^{2n}\label{Orderingdirectsum}\\
(e_{i},0)&\mapsto e_{2i-1}\nonumber\\
(0,e_{i})&\mapsto e_{2i}\nonumber.
\end{align}
The composition \eqref{OplusforOVFf} naturally induces a map 
\begin{equation}\label{OplusforMVFf}
\operatorname{Mor}^{m,n}_{k.l}(\mathbb{F})\times \operatorname{Mor}^{m,n}_{k^{\prime},l^{\prime}}(\mathbb{F})\rightarrow \operatorname{Mor}^{2m,2n}_{k+l,k^{\prime}+l^{\prime}}(\mathbb{F}), 
\end{equation} 
and the continuity of the maps \eqref{OplusforOVFf} and \eqref{OplusforMVFf} can be verified by examining their restriction to the coordinate charts. After properly ordering the bases, the restrictions of the maps \eqref{OplusforOVFf} and \eqref{OplusforMVFf} to the neighborhood $\mathfrak{N}_{X_{0}}\times \mathfrak{N}_{X^{\prime}_{0}}$ and $\mathfrak{N}_{X_{0},Y_{0}}\times \mathfrak{N}_{X_{0}^{\prime},Y_{0}^{\prime}}$ can be realized by the assignments
\[(A_{X},B_{X^{\prime}})\rightarrow 
\begin{bmatrix}
A_{X}&0\\
0&B_{X^{\prime}}
\end{bmatrix},
\] 
and  
\[(A_{X},A_{Y},[T]; B_{X^{\prime}},B_{Y^{\prime}},[S])\rightarrow (
\begin{bmatrix}
A_{X}&0\\
0&B_{X^{\prime}}
\end{bmatrix}
,
\begin{bmatrix}
A_{Y}&0\\
0&B_{Y^{\prime}}
\end{bmatrix}
,
\begin{bmatrix}
[T]&0\\
0&[S]
\end{bmatrix}
)
,\] 
respectively. Furthermore, from the way we order the basis of the ambient, meaning the assignment \eqref{Orderingdirectsum}, it is not difficult to see the maps \eqref{OplusforOVFf} and \eqref{OplusforMVFf} respect the stabilization maps, namely that the following diagrams commute 
\begin{center}
\begin{tikzpicture}
\node(Lu) at (0,2) {$\operatorname{Gr}_{k}^{n}\times \operatorname{Gr}_{k^{\prime}}^{n}$};
\node(Ll) at (0,0) {$\operatorname{Gr}_{k}^{n+1}\times \operatorname{Gr}_{k^{\prime}}^{n+1}$}; 
\node(Ru) at (5,2) {$\operatorname{Gr}_{k+k^{\prime}}^{2n}$};
\node(Rl) at (5,0) {$\operatorname{Gr}_{k+k^{\prime}}^{2(n+1)}$};

\path[->, font=\scriptsize,>=angle 90]

(Lu) edge (Ru)  
(Lu) edge (Ll)
(Ll) edge (Rl) 
(Ru) edge (Rl);

\end{tikzpicture}
\end{center}
and 
\begin{center}
\begin{tikzpicture}
\node(Lu) at (0,2) {$\operatorname{Mor}_{k,l}^{m,n}\times\operatorname{Mor}_{k^{\prime},l^{\prime}}^{m,n}$};
\node(Ll) at (0,0) {$\operatorname{Mor}_{k,l}^{m+1,n+1}\times\operatorname{Mor}_{k^{\prime},l^{\prime}}^{m+1,n+1}$}; 
\node(Ru) at (7,2) {$\operatorname{Mor}_{kk^{\prime},ll^{\prime}}^{2m,2n}$};
\node(Rl) at (7,0) {$\operatorname{Mor}_{kk^{\prime},ll^{\prime}}^{2(m+1),2(n+1)}$};

\path[->, font=\scriptsize,>=angle 90]

(Lu) edge (Ru)  
(Lu) edge (Ll)
(Ll) edge (Rl) 
(Ru) edge (Rl);

\end{tikzpicture}
\end{center}
where $\operatorname{Gr}_{k}^{n}$ and $\operatorname{Mor}_{k,l}^{m,n}$ denote the Grassmannian $\operatorname{Gr}_{k}(\mathbb{F}^{n})$ and the manifold $\operatorname{Mor}_{k,l}^{m,n}(\mathbb{F})$, respectively---for the sake of simplicity, we sometimes drop $\mathbb{F}$ from the notation. Passing to the colimit, we get a functor
\[\oplus:\mathcal{V}_{\mathbb{F}}^{f}\times  \mathcal{V}_{\mathbb{F}}^{f}\rightarrow \mathcal{V}_{\mathbb{F}}^{f}.\]

Consider now the tensor product of vector spaces and the  isomorphism given by  
\begin{align}\label{Orderingtensor}
\kappa:\mathbb{F}^{n}\otimes \mathbb{F}^{n}&\mapsto \mathbb{F}^{n^{2}}\nonumber\\
e_{i}\otimes e_{n}&\mapsto e_{i+(n-1)^{2}}  & i\leq n\\
e_{n}\otimes e_{j}&\mapsto e_{2n-j+(n-1)^{2}} & j\leq n.\nonumber
\end{align}
The following figure illustrates the ordering of the basis of $\mathbb{F}^{n^{2}}$.
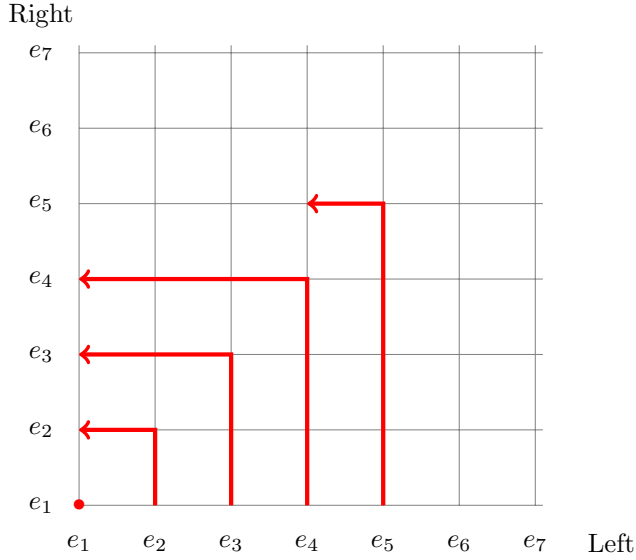
\begin{center}
\begin{tikzpicture}
\draw[gray,very thin] (1,1) grid (7.1,7.1);

\foreach \x in {1,2,3,4,5,6,7}
    \node at (\x, .5) {$e_{\x} $};
\foreach \y in {1,2,3,4,5,6,7}
    \node at (.5,\y) {$e_{\y}$};
    
\node[red] at (1,1){$\bullet$};    
\node at (8,.5) {$\operatorname{Left}$};
\node at (.5,7.5) {$\operatorname{Right}$};
\draw[red,ultra thick, ->] (2,1)--(2,2)--(1,2);
\draw[red,ultra thick, ->] (3,1)--(3,2)--(3,3)--(2,3)--(1,3);
\draw[red,ultra thick, ->] (4,1)--(4,2)--(4,3)--(4,4)--(3,4)--(2,4)--(1,4);
\draw[red,->,ultra thick, ->] (5,1)--(5,2)--(5,3)--(5,4)--(5,5)--(4,5);
    
\end{tikzpicture}
\end{center}
The isomorphism \eqref{Orderingtensor} induces a map 
\begin{equation}\label{OtimesOVFf}
\operatorname{Gr}_{k}^{n}\times \operatorname{Gr}_{k^{\prime}}^{n}\rightarrow \operatorname{Gr}_{kk^{\prime}}^{n^{2}}
\end{equation}
given by  
\begin{align*}
\mathbb{F}^{n}\times \mathbb{F}^{n}&\rightarrow \mathbb{F}^{n}\otimes\mathbb{F}^{n}\rightarrow \mathbb{F}^{n^{2}}\\
(V, W)&\mapsto  V\otimes W \mapsto \kappa_{\ast}(V\otimes W)  
\end{align*}
and a map 
\begin{align}\label{OtimesMVFf}
\operatorname{Mor}_{k,l}^{m,n}\times \operatorname{Mor}_{k^{\prime},l^{\prime}}^{m,n}&\rightarrow \operatorname{Mor}_{kk^{\prime},ll^{\prime}}^{m^{2},n^{2}}\\
(V\xrightarrow{T} W ; V^{\prime}\xrightarrow{T^{\prime}} W^{\prime})&\mapsto (\kappa_{\ast}(V\otimes V^{\prime})\xrightarrow{\kappa_{\ast}(T\otimes T^{\prime})} \kappa_{\ast}(W\otimes W^{\prime}))\nonumber.
\end{align} 
Similar to the case of the direct sum of vector spaces, one can check the continuity of the maps \eqref{OtimesOVFf} and \eqref{OtimesMVFf} by examining their restrictions to a neighborhood. For example, the map \eqref{OtimesOVFf}, when restricted to the neighborhood
$\mathfrak{N}_{X_{0}}\times \mathfrak{N}_{X^{\prime}_{0}}$, can be realized by  
\[(A_{X},A_{X^{\prime}})\mapsto 
A_{\kappa_{\ast}(X\otimes X^{\prime})}=
\left[
\begin{array}{c} 
A_{X}\otimes \operatorname{id}  \\
\hline
\operatorname{id}\otimes A_{X^{\prime}}\\
\hline
A_{X}\otimes A_{X^{\prime}}
\end{array}
\right],
\]
where the basis of $\mathbb{F}^{n^{2}}$ used for the neighborhood $\mathfrak{N}_{\kappa_{\ast}(X_{0}\otimes X^{\prime}_{0})}$ is   
\[ [\mathbf{x_{0}}\otimes \mathbf{x}_{0}^{\prime}, \mathbf{x}_{0}\otimes \mathbf{x}_{0}^{\prime,\perp},\mathbf{x}_{0}^{\perp}\otimes \mathbf{x}_{0}^{\prime},\mathbf{x}_{0}^{\perp}\otimes\mathbf{x}_{0}^{\prime,\perp}],\]
and $\mathbf{x}_{0}$ and $\mathbf{x}_{0}^{\prime}$ (resp. $\mathbf{x}_{0}^{\prime}$ and $\mathbf{x}_{0}^{\prime,\perp}$) are the bases for $X_{0}$ and $X_{0}^{\perp}$ (resp. $X_{0}^{\prime}$ and $X_{0}^{\prime,\perp}$), respectively. In the same way, one can check the continuity of the map \eqref{OtimesMVFf}
\[\operatorname{Mor}_{k,l}^{m,n}\times\operatorname{Mor}_{k^{\prime},l^{\prime}}^{m,n}\rightarrow\operatorname{Mor}_{kk^{\prime},ll^{\prime}}^{m^{2},n^{2}}\] by looking at its restriction on $\mathfrak{N}_{X_{0},Y_{0}}\times\mathfrak{N}_{X^{\prime}_{0},Y^{\prime}_{0}}$, which can be realized by
\begin{align*}
\mathfrak{N}_{X_{0},Y_{0}}\times\mathfrak{N}_{X^{\prime}_{0},Y^{\prime}_{0}}&\rightarrow\mathfrak{N}_{\kappa_{\ast}(X_{0}\otimes X^{\prime}_{0}),\kappa_{\ast}(Y_{0}\otimes Y^{\prime}_{0})},\\
\{(A_{X},B_{Y},[T]),(A_{X^{\prime}},B_{Y^{\prime}},[S])\}&\mapsto  (A_{\kappa_{\ast}(X\otimes X^{\prime})},B_{\kappa_{\ast}(Y\otimes Y^{\prime})},\kappa_{\ast}([T]\otimes [S])).
\end{align*} 
In view of the definition of the isomorphism $\mathbb{F}^{n}\otimes \mathbb{F}^{n}\simeq \mathbb{F}^{n^{2}}$ in \eqref{Orderingtensor}, it is not difficult to see they are compatible with the stablization maps, meaning the commutative diagrams
\begin{center}
\begin{tikzpicture}
\node(Lu) at (0,2) {$\operatorname{Gr}_{k}^{n}\times \operatorname{Gr}_{k^{\prime}}^{n}$};
\node(Ll) at (0,0) {$\operatorname{Gr}_{k}^{n+1}\times \operatorname{Gr}_{k^{\prime}}^{n+1}$}; 
\node(Ru) at (5,2) {$\operatorname{Gr}_{kk^{\prime}}^{n^{2}}$};
\node(Rl) at (5,0) {$\operatorname{Gr}_{kk^{\prime}}^{(n+1)^{2}}$};

\path[->, font=\scriptsize,>=angle 90]

(Lu) edge (Ru)  
(Lu) edge (Ll)
(Ll) edge (Rl) 
(Ru) edge (Rl);

\end{tikzpicture}
\end{center}
and 
\begin{center}
\begin{tikzpicture}
\node(Lu) at (0,2) {$\operatorname{Mor}_{k,l}^{m,n}\times\operatorname{Mor}_{k^{\prime},l^{\prime}}^{m,n}$};
\node(Ll) at (0,0) {$\operatorname{Mor}_{k,l}^{m+1,n+1}\times\operatorname{Mor}_{k^{\prime},l^{\prime}}^{m+1,n+1}$}; 
\node(Ru) at (7,2) {$\operatorname{Mor}_{kk^{\prime},ll^{\prime}}^{m^{2},n^{2}}$};
\node(Rl) at (7,0) {$\operatorname{Mor}_{kk^{\prime},ll^{\prime}}^{(m+1)^{2},(n+1)^{2}}$};

\path[->, font=\scriptsize,>=angle 90]

(Lu) edge (Ru)  
(Lu) edge (Ll)
(Ll) edge (Rl) 
(Ru) edge (Rl);

\end{tikzpicture}
\end{center}
Applying the colimit, we obtain a functor 
\[\otimes:\mathcal{V}_{\mathbb{F}}^{f}\times \mathcal{V}_{\mathbb{F}}^{f}\rightarrow \mathcal{V}_{\mathbb{F}}^{f}.\]
 
\subsection{H-semiring structure on $\vert iso\mathcal{V}_{F,\cdot}^{f} \vert$} 
Here we consider the internal subcategory of isomorphisms of $\mathcal{V}_{\mathbb{F}}^{f}$, whose space of objects is the same as that of $\mathcal{V}_{\mathbb{F}}^{f}$ and whose space of morphisms is the subspace of $M(\mathcal{V}_{\mathbb{F}}^{f})$ consisting of isomorphisms of vector spaces. It is denoted by  
$iso\mathcal{V}_{\mathbb{F}}^{f}$.  

Now, consider the subspace $iso_{k,k}^{m,n}(\mathbb{F})\subset \operatorname{Mor}_{k,k}^{m,n}(\mathbb{F})$ that comprises isomorphisms from a $k$-dimensional subspace in $\mathbb{F}^{m}$ to a $k$-dimensional subspace in $\mathbb{F}^{n}$, and observe that there is a bundle map
\begin{equation}\label{bundleofiso}
iso(s^{m}_{k},t^{n}_{k}):iso_{k,k}^{m,n}(\mathbb{F})\rightarrow \operatorname{Gr}_{k}(\mathbb{F}^{m})\times 
\operatorname{Gr}_{k}(\mathbb{F}^{n}),
\end{equation}
which can be viewed a subspace of $\operatorname{Mor}_{k,k}^{m,n}(\mathbb{F})$ over $\operatorname{Gr}_{k}(\mathbb{F}^{m})\times \operatorname{Gr}_{k}(\mathbb{F}^{n})$. In particular, at a point $(X_{0},Y_{0})$, the coordinate chart \eqref{CoordinateMor} of the vector bundle
\[(s^{m}_{k},t^{n}_{k}):\operatorname{Mor}_{k,k}^{m,n}(\mathbb{F})\rightarrow \operatorname{Gr}_{k}(\mathbb{F}^{m})\times \operatorname{Gr}_{k}(\mathbb{F}^{n})\]
gives us a diffeomorphism  
\[iso(s^{m}_{k},t^{n}_{k})^{-1}(\mathfrak{N}_{X_{0}}\times\mathfrak{N}_{Y_{0}})\rightarrow
M_{m-k,k}(\mathbb{F})\times M_{n-k,k}(\mathbb{F})\times \operatorname{GL}_{k}(\mathbb{F}).\]
Hence, the bundle \eqref{bundleofiso} is a fiber bundle with fiber diffeomorphic to $\operatorname{GL}_{k}(\mathbb{F})$.

\begin{lemma}\label{IsoVft}
The canonical map
\[\coprod_{k}\operatorname*{colim}_{m,n}iso_{k,k}^{m,n}(\mathbb{F})\rightarrow M(iso\mathcal{V}_{\mathbb{F}}^{f})\]
is a homeomorphism.
\end{lemma} 
\begin{proof}
The canonical map is clearly a bijection, so it suffices to show the subspace topology of $M(iso\mathcal{V}_{\mathbb{F}}^{f})$ in $\mathpzc{Top}^{w}$ coincides with the colimit topology induced by
$\{iso_{k,k}^{n,n}(\mathbb{F})\}_{n}$.  

To see this, we first observe that the space $\operatorname{Mor}^{n,n}_{k,k}(\mathbb{F})$ is normal, and secondly, we recall that the map 
\[\operatorname{Mor}^{n,n}_{k,k}(\mathbb{F})\rightarrow \operatorname{Mor}^{n+1,n+1}_{k,k}(\mathbb{F})\] 
is a cofibration, for every $n$. Furthermore, we have the intersection  
\[iso\mathcal{V}_{\mathbb{F}}^{f}\cap \operatorname{Mor}_{k,k}^{n,n}(\mathbb{F})=iso_{k,k}^{n,n}(\mathbb{F})\] 
in $\mathcal{T}$, the category of topological spaces, is a weak Hausdorff $k$-space. Therefore, the conditions in Lemma \ref{Colimandsubspace} are satisfied and the two topologies should coincide. 
\end{proof}

Now, the restrictions of the functors $\oplus$ and $\otimes$ to $iso\mathcal{V}_{\mathbb{F}}^{f}$ give us:
\begin{align*}
\oplus:iso\mathcal{V}_{\mathbb{F}}^{f}\times iso\mathcal{V}_{\mathbb{F}}^{f}&\rightarrow iso\mathcal{V}_{\mathbb{F}}^{f}\\
 \otimes:iso\mathcal{V}_{\mathbb{F}}^{f}\times iso\mathcal{V}_{\mathbb{F}}^{f}&\rightarrow iso\mathcal{V}_{\mathbb{F}}^{f}
\end{align*}
which, after applying the nerve construction and the realization functor, yield two operations on $\vert \operatorname{Ner}_{\cdot}iso\mathcal{V}_{\mathbb{F}}^{f}\vert$, denoted still by $\oplus$ and $\otimes$.  

\begin{theorem}\label{HsemiringVFf}
The space $\vert \operatorname{Ner}_{\cdot}iso\mathcal{V}_{\mathbb{F}}^{f}\vert$ is endowed with a $H$-semiring space structure with additive and multiplicative operations given by 
\begin{align*}
\oplus:\vert \operatorname{Ner}_{\cdot}iso\mathcal{V}_{\mathbb{F}}^{f}\vert\times \vert \operatorname{Ner}_{\cdot}iso\mathcal{V}_{\mathbb{F}}^{f}\vert &\rightarrow \vert \operatorname{Ner}_{\cdot}iso\mathcal{V}_{\mathbb{F}}^{f}\vert;\\
\otimes:\vert \operatorname{Ner}_{\cdot}iso\mathcal{V}_{\mathbb{F}}^{f}\vert \times \vert \operatorname{Ner}_{\cdot}iso\mathcal{V}_{\mathbb{F}}^{f}\vert &\rightarrow \vert\operatorname{Ner}_{\cdot} iso\mathcal{V}_{\mathbb{F}}^{f}\vert,
\end{align*}
respectively, and the additive and multiplicative identities by the points 
\begin{align*}
\ast &=\operatorname{Gr}_{0}(\mathbb{F}^{\infty});\\
\mathbb{F} =<e_{1}>&\in\operatorname{Gr}_{1}(\mathbb{F}^{\infty}),
\end{align*} 
respectively.
\end{theorem}
\begin{proof}
By the definition of $\oplus$, namely the isomorphism \eqref{Orderingdirectsum}, the diagrams 
\begin{center}
\begin{tikzpicture}
\node(Lu) at (0,2) {$\operatorname{Gr}_{i}^{n}\times\operatorname{Gr}_{j}^{n}\times \operatorname{Gr}_{k}^{n}$};
\node(Ll) at (0,0) {$\operatorname{Gr}_{i}^{n}\times \operatorname{Gr}_{j+k}^{2n}$}; 
\node(Ru) at (7,2) {$\operatorname{Gr}_{i+j}^{2n}\times \operatorname{Gr}_{k}^{n}$};
\node(Rl) at (7,0) {$\operatorname{Gr}_{i+j+k}^{3n}$};

\path[->, font=\scriptsize,>=angle 90]

(Lu) edge node [above]{$(\oplus,\operatorname{id})$}(Ru)  
(Lu) edge node [left]{$(\operatorname{id},\oplus)$}(Ll)
(Ll) edge node [above]{$\oplus$}(Rl) 
(Ru) edge node [right]{$\oplus$}(Rl);

\end{tikzpicture}
\end{center}
and
\begin{center}
\begin{tikzpicture}
\node(Lu) at (0,2) {$iso_{i,i}^{n,n}\times iso_{j,j}^{n,n}\times iso_{k,k}^{n,n}$};
\node(Ll) at (0,0) {$iso_{i,i}^{n,n}\times iso_{j+k,j+k}^{2n,2n}$}; 
\node(Ru) at (7,2) {$iso_{i+j}^{2n,2n} \times iso_{k,k}^{n,n} $};
\node(Rl) at (7,0) {$iso_{i+j+k,i+j+k}^{3n,3n} $};

\path[->, font=\scriptsize,>=angle 90]

(Lu) edge (Ru)  
(Lu) edge (Ll)
(Ll) edge (Rl) 
(Ru) edge (Rl);

\end{tikzpicture}
\end{center}
are commutative, and hence the associativity of $\oplus$ follows from Lemma \ref{Colimitandpullback}. Note that, in the diagrams above, $\mathbb{F}$ is dropped from the notation for the sake of simplicity.

To check the point $\ast$ is indeed the additive identity, we observe that the following diagram
\begin{center}
\begin{equation}\label{AddIdinVFf}
\begin{tikzpicture}[baseline=(current bounding box.center)]
\node(Lu) at (0,2) {$\operatorname{Gr}_{k}^{n}$};
\node(Ll) at (0,0) {$\operatorname{Gr}_{0}^{n}\times \operatorname{Gr}_{k}^{n}$}; 
\node(Rl) at (5,0) {$\operatorname{Gr}_{k}^{2n}$};

\path[->, font=\scriptsize,>=angle 90]

(Lu) edge node [left]{$(\operatorname{id},\ast)/(\ast,\operatorname{id})$}(Ll)
(Ll) edge node [above]{$\oplus$}(Rl) 
(Lu) edge node [above]{$\operatorname{id}$}(Rl);

\end{tikzpicture}
\end{equation}
\end{center} 
commutes up to natural maps. More precisely, if we denote the image of $X\in \operatorname{Gr}_{k}(\mathbb{F}^{n})$ under the map 
\[\operatorname{Gr}_{k}(\mathbb{F}^{n})\rightarrow\operatorname{Gr}_{k}(\mathbb{F}^{2n})\]
by $X^{\prime}$ and the image of $X$ under the composition 
\[\operatorname{Gr}_{k}(\mathbb{F}^{n})\rightarrow\operatorname{Gr}_{0}(\mathbb{F}^{n})\times\operatorname{Gr}_{k}(\mathbb{F}^{n})\xrightarrow{\oplus}\operatorname{Gr}_{k}(\mathbb{F}^{2n})\]
by $X^{\prime\prime}$, then the identity of $X$ gives a natural map from $X^{\prime}$ to $X^{\prime\prime}$ in $\operatorname{Mor}_{k,k}^{2n,2n}(\mathbb{F})$. If one orders the bases of the subspaces $X^{\prime}_{0}$, $X^{\prime\prime}_{0}$, $X^{\prime,\perp}_{0}$ and $X^{\prime\prime,\perp}_{0}$ in $\mathbb{F}^{2n}$ properly, the natural map can be realized in terms of the coordinate charts by the assignment  
\begin{align*}
\mathfrak{N}_{X_{0}}\simeq M_{n-k,k}(\mathbb{F})& \rightarrow M_{2n-k,k}(\mathbb{F})\times M_{2n-k,k}(\mathbb{F})\times M_{k,k}(\mathbb{F})\times \operatorname{GL}_{k}(\mathbb{F})\simeq \mathfrak{N}_{X_{0}^{\prime},X^{\prime\prime}_{0}}\\
A_{X}&\mapsto ([A_{X},0],[A_{X},0],[\operatorname{id}])
\end{align*} 
By passing to the colimit and Lemma \ref{Colimitandpullback}, we see the diagram  
\begin{center}
\begin{tikzpicture}
\node(Lu) at (0,2) {$iso\mathcal{V}_{\mathbb{F}}^{f}$};
\node(Ll) at (0,0) {$iso\mathcal{V}_{\mathbb{F}}^{f}\times iso\mathcal{V}_{\mathbb{F}}^{f}$}; 
\node(Rl) at (5,0) {$iso\mathcal{V}_{\mathbb{F}}^{f}$};

\path[->, font=\scriptsize,>=angle 90]

(Lu) edge node [left]{$(\operatorname{id},\ast)/(\ast,\operatorname{id})$}(Ll)
(Ll) edge node [above]{$\oplus$}(Rl) 
(Lu) edge node [above]{$\operatorname{id}$}(Rl);

\end{tikzpicture}
\end{center} 
commutes up to natural transformations, and hence the diagram
\begin{center}
\begin{tikzpicture}
\node(Lu) at (0,2) {$\vert\operatorname{Ner}_{\cdot}iso\mathcal{V}_{\mathbb{F}}^{f}\vert$};
\node(Ll) at (0,0) {$\vert\operatorname{Ner}_{\cdot}iso\mathcal{V}_{\mathbb{F}}^{f}\vert \times \vert\operatorname{Ner}_{\cdot}iso\mathcal{V}_{\mathbb{F}}^{f}\vert$}; 
\node(Rl) at (5,0) {$\vert\operatorname{Ner}_{\cdot}iso\mathcal{V}_{\mathbb{F}}^{f}\vert$};

\path[->, font=\scriptsize,>=angle 90]

(Lu) edge node [left]{$(\operatorname{id},\ast)/(\ast,\operatorname{id})$}(Ll)
(Ll) edge node [above]{$\oplus$}(Rl) 
(Lu) edge node [above]{$\operatorname{id}$}(Rl);
\end{tikzpicture}
\end{center} 
commutes up to homotopy, in view of Construction \ref{ConstinIntCat}.
 
Similarly, we can easily verify that the diagram  
\begin{center}
\begin{tikzpicture}
\node(Lu) at (0,2) {$\operatorname{Gr}_{k}^{n}\times \operatorname{Gr}_{k}^{n}$};
\node(Ml) at (2,0) {$\operatorname{Gr}_{k}^{2n}$}; 
\node(Ru) at (4,2) {$\operatorname{Gr}_{k}^{n}\times \operatorname{Gr}_{k}^{n}$};
 
\path[->, font=\scriptsize,>=angle 90]

(Lu) edge node [above]{$\tau$}(Ru)  
(Lu) edge node [left]{$\oplus$}(Ml)
(Ru) edge node [right]{$\oplus$}(Ml);

\end{tikzpicture}
\end{center}
commutes up to natural maps, where $\tau$ sends $(x,y)$ to $(y,x)$ and the natural map from $\theta_{\ast}(X\oplus Y)$ to $\theta_{\ast}(Y\oplus X)$ is induced by the natural isomorphism $X\oplus Y\simeq Y\oplus X$. If one choose and order the bases of $\mathbb{F}^{2n}$ properly, then the natural map can be realized by the assignment below, in terms of the coordinate charts.
\begin{align*} 
M_{n-k,k}(\mathbb{F})\times M_{n-k,k}(\mathbb{F}) \rightarrow& M_{2n-2k,2k}(\mathbb{F})\times M_{2n-2k,2k}(\mathbb{F})\times\operatorname{GL}_{2k}(\mathbb{F})\\
(A_{X},A_{X^{\prime}})\mapsto& (\theta_{\ast}(A_{X}\oplus A_{X^{\prime}}),\theta_{\ast}(A_{X^{\prime}}\oplus A_{X}),[\frac{0,\operatorname{id}}{\operatorname{id},0}]).
\end{align*}
This shows the diagram    
\begin{center}
\begin{tikzpicture}
\node(Lu) at (0,2) {$iso\mathcal{V}_{\mathbb{F}}^{f}\times iso\mathcal{V}_{\mathbb{F}}^{f}$};
\node(Ml) at (2,0) {$iso\mathcal{V}_{\mathbb{F}}^{f}$}; 
\node(Ru) at (4,2) {$iso\mathcal{V}_{\mathbb{F}}^{f}\times iso\mathcal{V}_{\mathbb{F}}^{f}$};
 
\path[->, font=\scriptsize,>=angle 90]

(Lu) edge node [above]{$\tau$}(Ru)  
(Lu) edge node [left]{$\oplus$}(Ml)
(Ru) edge node [right]{$\oplus$}(Ml);

\end{tikzpicture}
\end{center}
commutes up to natural transformations and hence the operation $\oplus$ is commutative.

The proofs of the commutativity and associativity of the multiplicative operation $\oplus$ and the point $\mathbb{F}$ is the multiplication identity can be worked out in exactly the same way as that of the additive operation. Namely, one first checks the the diagram on the finite level commutes up to natural maps, and then passing to the colimit to get a commutative diagram up to natural transformations. One thing to note is that the operation $\otimes$ descends to a map from the smash product
\[\otimes:\vert iso\mathcal{V}_{\mathbb{F}}^{f}\vert\wedge \vert  iso\mathcal{V}_{\mathbb{F}}^{f}\vert\rightarrow \vert iso\mathcal{V}_{\mathbb{F}}^{f}\vert.\]

In the same vein, for the distribution law, the natural isomorphisms 
\begin{align*}
(X\oplus Y)\otimes Z\simeq (X\otimes Z)\oplus (Y\otimes Z),\\ 
X\otimes (Y\oplus Z)\simeq (X\otimes Y)\oplus (X\otimes Z)
\end{align*}
induce the natural transformations needed for the diagram 
\begin{center}
\begin{tikzpicture}
\node(Lu) at (0,2) {$iso\mathcal{V}_{\mathbb{F}}^{f}\times iso\mathcal{V}_{\mathbb{F}}^{f}\times iso\mathcal{V}_{\mathbb{F}}^{f}$};
\node(Lm) at (0,1) {$iso\mathcal{V}_{\mathbb{F}}^{f}\times iso\mathcal{V}_{\mathbb{F}}^{f}\times iso\mathcal{V}_{\mathbb{F}}^{f}\times iso\mathcal{V}_{\mathbb{F}}^{f}$};
\node(Ll) at (0,0) {$iso\mathcal{V}_{\mathbb{F}}^{f}\times iso\mathcal{V}_{\mathbb{F}}^{f}$}; 
\node(Ru) at (5,2) {$iso\mathcal{V}_{\mathbb{F}}^{f}\times iso\mathcal{V}_{\mathbb{F}}^{f}$};
 
\node(Rl) at (5,0) {$iso\mathcal{V}_{\mathbb{F}}^{f}$};

\path[->, font=\scriptsize,>=angle 90]

(Lu) edge node [above]{$(\operatorname{id},\oplus)/(\oplus,\operatorname{id})$}(Ru)  
(Lu) edge node [right]{$l/r$}(Lm)
(Lm) edge node [right]{$(\otimes,\otimes)$}(Ll)
(Ll) edge node [above]{$\oplus$}(Rl) 
(Ru) edge node [right]{$\otimes$}(Rl);

\end{tikzpicture}
\end{center} 
where $l(x,y,z)=(x,y,x,z)$ and $r(x,y,z)=(x,z,y,z)$.  

\end{proof}

\subsection{Geometric meaning of $\mathcal{V}_{\mathbb{F}}^{f}$}
Here, we shall examine some geometric and topological properties of the category $iso\mathcal{V}_{\mathbb{F}}^{f}$. In particular, we shall see the space
\[\vert\operatorname{Ner}_{\cdot}iso\mathcal{V}_{\mathbb{F}}^{f}\vert\]
classifies isomorphism classes of vector bundles.

\begin{lemma}\label{Universalbundleofisom}
Let \[\operatorname{E}_{k} \rightarrow\operatorname{Gr}_{k}(\mathbb{F}^{\infty})\] 
be the universal $k$-dimensional vector bundle and $iso_{k}(\mathbb{F})$ denote the colimit 
\[\operatorname*{colim}_{n} iso_{k,k}^{n,n}(\mathbb{F}).\]
Then, given a map  
\[X\xrightarrow{(f,g)} \operatorname{Gr}_{k}(\mathbb{F}^{\infty})\times \operatorname{Gr}_{k}(\mathbb{F}^{\infty}),\]
there is a bundle isomorphism 
\begin{center}
\begin{equation}\label{MorclassifiesBundleiso}
\begin{tikzpicture}[baseline=(current bounding box.center)] 
\node(Lu) at (0,2) {$(f,g)^{\ast}iso_{k}(\mathbb{F})$};
\node(Ru) at (6,2) {$\operatorname{Isom}(f^{\ast}\operatorname{E}_{k},g^{\ast}\operatorname{E}_{k})$};
\node(Ml) at (3,0) {$X$};

\path[->, font=\scriptsize,>=angle 90]

(Lu) edge (Ru)  
(Lu) edge (Ml)  
(Ru) edge (Ml);

\end{tikzpicture}
\end{equation}
\end{center}
where $(f,g)^{\ast}iso_{k}(\mathbb{F})$ is the pullback of the principal $\operatorname{GL}_{k}(\mathbb{F})$-bundle 
\[iso_{k}(\mathbb{F})\rightarrow \operatorname{Gr}_{k}(\mathbb{F}^{\infty})\times \operatorname{Gr}_{k}(\mathbb{F}^{\infty})\]
along $(f,g)$ and $\operatorname{Isom}(f^{\ast}\operatorname{E}_{k},g^{\ast}\operatorname{E}_{k})$ is the principal $\operatorname{GL}_{k}$-bundle of isomorphisms between $f^{\ast}\operatorname{E}_{k}$ and $g^{\ast}\operatorname{E}_{k}$.
\end{lemma}
\begin{proof} 
By definition, the diagram \eqref{MorclassifiesBundleiso} is an isomorphism over $X$ as sets. To see their topologies also coincide, we examine the coordinate charts. The coordinate charts of $(f,g)^{\ast}iso_{k}(\mathbb{F})$ are induced from the coordinate charts of $ iso_{k}(\mathbb{F})$, whereas the coordinate charts of $\operatorname{Isom}(f^{\ast}\operatorname{E}_{k},g^{\ast}\operatorname{E}_{k})$ are given by the coordinate charts of $f^{\ast}\operatorname{E}_{k}$ and $g^{\ast}\operatorname{E}_{k}$. From the construction \eqref{CoordinateMor}, we see they are exactly the same one.
\end{proof}

\begin{lemma}\label{Manyrealizations}
The canonical maps  
\[Biso\mathcal{V}_{\mathbb{F}}^{f}\rightarrow  \vert\vert \operatorname{Ner}_{\cdot}iso\mathcal{V}_{\mathbb{F}}^{f}\vert\vert\rightarrow \vert  \operatorname{Ner}_{\cdot}iso\mathcal{V}_{\mathbb{F}}^{f}\vert  \]
are homotopy equivalences, where the space $Biso\mathcal{V}_{\mathbb{F}}^{f}$ is the classifying of the topological groupoid $iso\mathcal{V}_{\mathbb{F}}^{f}$ constructed in \cite[Section $3$]{Se1} and \cite[Appendix $B$]{Bo} (see also \cite[Definition $3.1$]{Wang4}).
\end{lemma}
\begin{proof} 
By the variant of tom Dieck's theorem in \cite[Theorem $4.3$]{Wang4}, we know if the space $\operatorname{Ner}_{i}\mathcal{V}_{\mathbb{F}}^{f}$ has the homotopy type of a $\operatorname{CW}$-complex, then the map  
\begin{equation}\label{Map1inManyrealizations}
Biso\mathcal{V}_{\mathbb{F}}^{f}\rightarrow  \vert\vert \operatorname{Ner}_{\cdot}iso\mathcal{V}_{\mathbb{F}}^{f}\vert\vert
\end{equation} 
is a homotopy equivalence.   

Now, since the map
\[iso(s^{n}_{k},t^{n}_{k}):iso_{k,k}^{n,n}(\mathbb{F})\rightarrow \operatorname{Gr}_{k}(\mathbb{F}^{n})\times \operatorname{Gr}_{l}(\mathbb{F}^{n})\]
is a principle $\operatorname{GL}_{k}(\mathbb{F})$-bundle, we have the following pullback diagram
\begin{center}
\begin{tikzpicture}
\node (Lu) at (0,2){$iso_{k,k}^{n,n}(\mathbb{F})$};
\node (Ll) at (0,0){$\operatorname{Gr}_{k}(\mathbb{F}^{n})\times \operatorname{Gr}_{k}(\mathbb{F}^{n})$};
\node (Ru) at (4,2){$\operatorname{EGL}_{k}(\mathbb{F})$};
\node (Rl) at (4,0){$\operatorname{BGL}_{k}(\mathbb{F})$};

\path[->, font=\scriptsize,>=angle 90] 

(Lu) edge node[above] {$\hat{f}^{n}$}(Ru)
(Ll) edge node[above] {$f^{n}$}(Rl)
(Ru) edge node[right] {$\pi$}(Rl)
(Lu) edge node[right] {$iso(s^{n}_{k},t^{n}_{k})$}(Ll);

\end{tikzpicture}
\end{center}
where $\operatorname{EGL}_{k}(\mathbb{F})\rightarrow \operatorname{BGL}_{k}(\mathbb{F})$ is the universal $\operatorname{GL}_{k}(\mathbb{F})$ bundle constructed by the standard bar construction \cite[Definition 1.62]{Ru}, \cite{Ma1}, \cite{BV}. Since pullbacks commute with sequential colimits of injective maps (Lemma \ref{Colimitandpullback}), the diagram below 
\begin{center}
\begin{tikzpicture}
\node (Lu) at (0,2){$\operatorname*{colim}\limits_{n}iso_{k,k}^{n,n}(\mathbb{F})$};
\node (Ll) at (0,0){$\operatorname*{colim}\limits_{n}(\operatorname{Gr}_{k}(\mathbb{F}^{n})\times \operatorname{Gr}_{k}(\mathbb{F}^{n}))$};
\node (Ru) at (4,2){$\operatorname{EGL}_{k}(\mathbb{F})$};
\node (Rl) at (4,0){$\operatorname{BGL}_{k}$};

\path[->, font=\scriptsize,>=angle 90] 

(Lu) edge node[above] {$\operatorname*{colim}\limits_{n} \hat{f}^{n}$}(Ru)
(Ll) edge node[above] {$\operatorname*{colim}\limits_{n} f_{n}^{n}$}(Rl)
(Ru) edge node[right] {$\pi$}(Rl)
(Lu) edge node[right] {$\operatorname*{colim}\limits_{n} iso(s^{n}_{k},t^{n}_{k}) $}(Ll);

\end{tikzpicture}
\end{center} 
is also a pullback. In particular, it implies the map
\[M(iso\mathcal{V}_{\mathbb{F}}^{f})\rightarrow O(iso\mathcal{V}_{\mathbb{F}}^{f})\times O(iso\mathcal{V}_{\mathbb{F}}^{f})\] 
is a fiber bundle. Since both base space and fiber are $\operatorname{CW}$-complexes, the total space $M(iso\mathcal{V}_{\mathbb{F}}^{f})$ must have the homotopy type of a $\operatorname{CW}$-complex by Lemma \ref{TCWfibration}.
 
On the other hand, by the definition of the nerve construction, we have the pullback diagram
\begin{center}
\begin{tikzpicture}
\node (Lu) at (0,2){$\operatorname{Ner}_{i}iso\mathcal{V}_{\mathbb{F}}^{f}$};
\node (Ll) at (0,0){$\operatorname{Ner}_{i-1}iso\mathcal{V}_{\mathbb{F}}^{f}$};
\node (Ru) at (3,2){$M(iso\mathcal{V}_{\mathbb{F}}^{f})$};
\node (Rl) at (3,0){$O(iso\mathcal{V}_{\mathbb{F}}^{f})$};

\path[->, font=\scriptsize,>=angle 90] 

(Lu) edge node[above] {$p_{k}$}(Ru)
(Ll) edge node[above] {$t$}(Rl)
(Ru) edge node[right] {$s$}(Rl)
(Lu) edge node[right] {$p_{1,...,i-1}$}(Ll);

\end{tikzpicture}
\end{center} 
which implies the map
\begin{equation}\label{FiberationNerVFt}
\operatorname{Ner}_{i}iso\mathcal{V}_{\mathbb{F}}^{f}\rightarrow \operatorname{Ner}_{i-1}iso\mathcal{V}_{\mathbb{F}}^{f}
\end{equation} 
is also a Hurewicz fibration with fiber homeomorphic to the total space of a fiber bundle whose fiber and base space are the disjoint unions 
\[\coprod_{k}\operatorname{GL}_{k}(\mathbb{F})\text{ and }\coprod_{k}\operatorname{Gr}_{k}(\mathbb{F}^{\infty}),\]
respectively. Therefore, the fiber of the Hurewicz fibration \eqref{FiberationNerVFt} has the the homotopy type of a $\operatorname{CW}$-complex. By induction and Lemma \ref{TCWfibration}, we see that the space $\operatorname{Ner}_{i}iso\mathcal{V}_{\mathbb{F}}^{f}$ has the homotopy type of a $\operatorname{CW}$-complex, for every $i$. Hence, by Theorem $4.3$ in \cite{Wang4}, the map \eqref{Map1inManyrealizations} is a homotopy equivalence. 
 
Now, to show the map   
\[\vert\vert \operatorname{Ner}_{\cdot}iso\mathcal{V}_{\mathbb{F}}^{f}\vert\vert\rightarrow \vert  \operatorname{Ner}_{\cdot}iso\mathcal{V}_{\mathbb{F}}^{f}\vert\]
is a homotopy equivalence, we know, in view of Proposition $A.1.v$ in \cite{Se3}, it suffices to show that the simplicial space $\operatorname{Ner}_{\cdot}iso\mathcal{V}_{\mathbb{F}}^{f}$ is proper. To see this, we observe that the map 
\[e^{m}_{k}:\operatorname{Gr}_{k}(\mathbb{F}^{m})\rightarrow iso_{k,k}^{m,m}(\mathbb{F})\]
is a closed embedding of a smooth manifold and hence a cofibration. On the other hand, since the map
\begin{align*}
p_{r}\hspace*{.5em}(\operatorname{resp}. p_{0}):\overbrace{iso^{m,m}_{k,k}\times_{\operatorname{Gr}_{k}^{m}}...\times_{\operatorname{Gr}_{k}^{m}}iso_{k,k}^{m,m}}^{r}&\rightarrow \operatorname{Gr}_{k}^{m}\\
(V_{0}\xrightarrow{T_{1}}... \xrightarrow{T_{r}}V_{r})&\mapsto V_{r}\hspace*{.5em}(\operatorname{resp}. V_{0}) 
\end{align*}
is a composition of Hurewicz fibrations, it is itself a Hurewicz fibration, for any $r$. These two observations tell us the diagram below satisfy the conditions in Lemma \ref{PullbackandCof2}.  
\begin{center}
\begin{tikzpicture}
\node(Lu) at (0,2.2) {$\overbrace{iso^{m,m}_{k,k}\times_{\operatorname{Gr}_{k}^{m}}...\times_{\operatorname{Gr}_{k}^{m}}iso_{k,k}^{m,m}}^{r-1}$};
\node(Ll) at (0,0.2) {$\overbrace{iso^{m,m}_{k,k}\times_{\operatorname{Gr}_{k}^{m}}...\times_{\operatorname{Gr}_{k}^{m}}iso_{k,k}^{m,m}}^{r-1}$};
\node(Mu) at (4,2) {$\operatorname{Gr}_{k}^{m}$};
\node(Ml) at (4,0) {$\operatorname{Gr}_{k}^{m}$};
\node(Ru) at (6,2) {$\operatorname{Gr}_{k}^{m}$};
\node(Rl) at (6,0) {$iso^{m,m}_{k,k}$};

\path[->, font=\scriptsize,>=angle 90]

(Rl) edge node [above]{$s_{k}^{m}$}(Ml)
 
(Ru) edge node [above]{$=$}(Mu)  
(Lu) edge node [right]{$\parallel$}(Ll)
(Mu) edge node [right]{$\parallel$}(Ml)  
(Ru) edge node [right]{$e_{k}^{m}$}(Rl);
\draw[->](2.5,2) to node [above]{\scriptsize $p_{r-1}$}(Mu); 
\draw[->](2.5,0) to node [above]{\scriptsize $p_{r-1}$}(Ml);
\end{tikzpicture}
\end{center}
Thus, the map 
\begin{align*}
\mathbcal{s}_{r}^{m}:\overbrace{iso^{m,m}_{k,k}\times_{\operatorname{Gr}_{k}^{m}}...\times_{\operatorname{Gr}_{k}^{m}}iso_{k,k}^{m,m}}^{r}&\rightarrow \overbrace{iso^{m,m}_{k,k}\times_{\operatorname{Gr}_{k}^{m}}...\times_{\operatorname{Gr}_{k}^{m}}iso_{k,k}^{m,m}}^{r+1}\\
(T_{1},...,T_{r})&\mapsto (T_{1},...,T_{r},\operatorname{id}) 
\end{align*}
is a cofibration. It implies that the following map of cospans  
\begin{center}
\begin{tikzpicture}
\node(Lu) at (0,2.2) {$\overbrace{iso^{m,m}_{k,k}\times_{\operatorname{Gr}_{k}^{m}}...\times_{\operatorname{Gr}_{k}^{m}}iso_{k,k}^{m,m}}^{r^{\prime}-1}$};
\node(Ll) at (0,0.2) {$\overbrace{iso^{m,m}_{k,k}\times_{\operatorname{Gr}_{k}^{m}}...\times_{\operatorname{Gr}_{k}^{m}}iso_{k,k}^{m,m}}^{r^{\prime}}$};
\node(Mu) at (3.5,2) {$\operatorname{Gr}_{k}^{m}$};
\node(Ml) at (3.5,0) {$\operatorname{Gr}_{k}^{m}$};
\node(Ru) at (7,2.2) {$\overbrace{iso^{m,m}_{k,k}\times_{\operatorname{Gr}_{k}^{m}}...\times_{\operatorname{Gr}_{k}^{m}}iso_{k,k}^{m,m}}^{r^{\prime\prime}}$};
\node(Rl) at (7,0.2) {$\overbrace{iso^{m,m}_{k,k}\times_{\operatorname{Gr}_{k}^{m}}...\times_{\operatorname{Gr}_{k}^{m}}iso_{k,k}^{m,m}}^{r^{\prime\prime}}$};

\path[->, font=\scriptsize,>=angle 90]

(Lu) edge node [right]{$\mathbcal{s}^{m}_{r^{\prime}-1}$}(Ll)
(Mu) edge node [right]{$\parallel$}(Ml)  
(Ru) edge node [right]{$\parallel$}(Rl);

\draw[->](2.5,0) to node [above]{\scriptsize $p_{r^{\prime}}$}(Ml);
\draw[->](2.5,2) to node [above]{\scriptsize $p_{r^{\prime}-1}$}(Mu);

\draw[->](4.5,0) to node [above]{\scriptsize $p_{0}$}(Ml);
\draw[->](4.5,2) to node [above]{\scriptsize $p_{0}$}(Mu);
\end{tikzpicture}
\end{center}
also satisfies the conditions in Lemma \ref{PullbackandCof2}, and therefore, for every $i$, the map
\begin{align}\label{Cof1inManyrealization}
\mathbcal{s}_{i}^{m}:\overbrace{iso^{m,m}_{k,k}\times_{\operatorname{Gr}_{k}^{m}}...\times_{\operatorname{Gr}_{k}^{m}}iso_{k,k}^{m,m}}^{r}&\rightarrow \overbrace{iso^{m,m}_{k,k}\times_{\operatorname{Gr}_{k}^{m}}...\times_{\operatorname{Gr}_{k}^{m}}iso_{k,k}^{m,m}}^{r+1}\\
(T_{1},...,T_{r})&\mapsto (T_{1},..\overbrace{\operatorname{id}}^{i}.,T_{r})\nonumber 
\end{align}
is a cofibration. 

On the other hand, we know, for every $k,m$ and $r$, the map
\begin{multline}\label{Cof2inManyrealization}
\overbrace{iso^{m,m}_{k,k}\times_{\operatorname{Gr}_{k}^{m}}...\times_{\operatorname{Gr}_{k}^{m}}iso_{k,k}^{m,m}}^{r}\\
\rightarrow \overbrace{iso^{m+1,m+1}_{k,k}\times_{\operatorname{Gr}_{k}^{m+1}}...\times_{\operatorname{Gr}_{k}^{m+1}}iso_{k,k}^{m+1,m+1}}^{r}
\end{multline} 
is a cofibration---this can be seen by applying Lemma \ref{PullbackandCof2} to the diagram below and the induction
\begin{center}
\begin{tikzpicture}
\node(Lu) at (0,2.2) {$\overbrace{iso^{m,m}_{k,k}\times_{\operatorname{Gr}_{k}^{m}}...\times_{\operatorname{Gr}_{k}^{m}}iso_{k,k}^{m,m}}^{r-1}$};
\node(Ll) at (0,0.2) {$\overbrace{iso^{m+1,m+1}_{k,k}\times_{\operatorname{Gr}_{k}^{m+1}}...\times_{\operatorname{Gr}_{k}^{m+1}}iso_{k,k}^{m+1,m+1}}^{r-1}$};
\node(Mu) at (4.5,2) {$\operatorname{Gr}_{k}^{m}$};
\node(Ml) at (4.5,0) {$\operatorname{Gr}_{k}^{m+1}$};
\node(Ru) at (7,2) {$iso^{m,m}_{k,k}$};
\node(Rl) at (7,0) {$iso^{m+1,m+1}_{k,k}$};

\path[->, font=\scriptsize,>=angle 90]

(Rl) edge  (Ml) 
(Ru) edge  (Mu)  
(Lu) edge  (Ll)
(Mu) edge  (Ml)  
(Ru) edge  (Rl);
\draw[->](3.5,0) to  (Ml);
\draw[->](2.5,2) to  (Mu);

\end{tikzpicture}
\end{center}
Cofibrations \eqref{Cof1inManyrealization} and \eqref{Cof2inManyrealization} imply that the following map of sequences of maps
\[\{\mathbcal{s}_{i}^{m}\}_{m}:\{\overbrace{iso^{m,m}_{k,k}\times_{\operatorname{Gr}_{k}^{m}}...\times_{\operatorname{Gr}_{k}^{m}}iso_{k,k}^{m,m}}^{r-1}\}_{m}\mapsto \{\overbrace{iso^{m,m}_{k,k}\times_{\operatorname{Gr}_{k}^{m}}...\times_{\operatorname{Gr}_{k}^{m}}iso_{k,k}^{m,m}}^{r}\}_{m}\]
is a cofibrant object in the Reedy Model structure on the category of functors from $\mathbb{N}$ to $\mathpzc{Top}^{w}$. Hence, on passing to the colimit, we get all degeneracy maps in $\operatorname{Ner}_{\cdot}\mathcal{V}_{\mathbb{F}}^{f}$ are cofibrations.
\end{proof}

We are now in the position to prove the main theorem of this section.  
\begin{theorem}\label{bijectionhomotopyandvectorbundle}
For any paracompact Hausdorff space $X$, there is a $1$-$1$ correspondence:
\[ [X,\vert iso\mathcal{V}_{\mathbb{F}}^{f}\vert]\leftrightarrow \coprod_{k}\operatorname{Vect}^{\mathbb{F}}_{k}(X).\] 
\end{theorem}

\begin{proof} 
By Theorem $3.2$ in \cite{Wang4} and Lemma \ref{Manyrealizations}, we know that there is a bijection
\[[X,\vert iso\mathcal{V}_{\mathbb{F}}^{f}\vert]\leftrightarrow  iso\mathcal{V}_{\mathbb{F}}^{f}(X).\]
Hence, it suffices to show there is a $1$-$1$ correspondence between 
\[iso\mathcal{V}_{\mathbb{F}}^{f}(X)\leftrightarrow \coprod_{k}\operatorname{Vect}_{k}(X).\]
Without loss of generality, we may assume $X$ is connected. Given an $iso\mathcal{V}_{\mathbb{F}}^{f}$-cocycle $\{\mathfrak{U}_{\alpha},f_{\beta\alpha}\}$, the map  
\[f_{\alpha\alpha}:\mathfrak{U}_{\alpha}\rightarrow  \operatorname{Gr}_{k}(\mathbb{F}^{\infty})\] 
determines a vector bundle over $\mathfrak{U}_{\alpha}$, denoted by $U_{\alpha}$, and the map
\[f_{\beta\alpha}:\mathfrak{U}_{\alpha}\cap \mathfrak{U}_{\beta}\rightarrow iso_{k}(\mathbb{F})\]
determines a bundle morphism $F_{\beta\alpha}:U_{\alpha}\vert_{\mathfrak{U}_{\alpha}\cap \mathfrak{U}_{\beta}}\rightarrow U_{\beta}\vert_{\mathfrak{U}_{\alpha}\cap \mathfrak{U}_{\beta}}$. The cocycle condition on $\{f_{\beta\alpha}\}$ implies $F_{\gamma\beta}\circ F_{\beta\alpha}=F_{\gamma\alpha}$. The associated (numerable) vector bundle is given by gluing $\{U_{\alpha}\}$ together by $\{F_{\beta\alpha}\}$, namely the quotient space  
\[U:=\coprod_{\alpha}U_{\alpha}/\sim,\]  
where the equivalence $\sim$ identifies $v\in U_{\alpha}$ with $F_{\beta\alpha}(v)\in U_{\beta}$.

Given two cocycles $u$ $v$ on $X$ which are connected by another cocycle $w$ on $X\times I$, then by the construction just described, there is a vector bundle $W$ over $X\times I$ which restricts to $U$ and $V$, the vector bundles induced by $u$ and $v$, respectively. Since $X$ is paracompact Hausdorff, $X\times I$ is also paracompact Hausdorff and thus $W$ is numerable. It implies $U$ and $V$ are isomorphic and hence there is a well-defined function 
\[\Psi: iso\mathcal{V}_{\mathbb{F}}^{f}(X)\rightarrow \coprod_{k}\operatorname{Vect}_{k}(X).\] 
The surjectivity of $\Psi$ follows from the commutative diagram:
\begin{center}
\begin{tikzpicture}
\node(Lu) at (0,2) {$iso\mathcal{V}_{\mathbb{F}}^{f}(X)$};
\node(Ll) at (0,0) {$iso\mathcal{V}_{\mathbb{F}}(X)$}; 
\node(Ru) at (5,2) {$\coprod\limits_{k}\operatorname{Vect}_{k}(X)$};
 
\path[->, font=\scriptsize,>=angle 90]

(Lu) edge node [above]{$\Psi$}(Ru)  
(Ll) edge node [right]{$\iota$}(Lu);

\draw[transform canvas={yshift=0.5ex},->] (Ll) --(Ru) node[above,midway] {\tiny $\psi$};
\draw[transform canvas={yshift=-0.5ex},->](Ru) -- (Ll) node[below,midway] {\tiny $\psi^{-1}$}; 
\end{tikzpicture}
\end{center}
where $iso\mathcal{V}_{\mathbb{F}}\rightarrow iso\mathcal{V}_{\mathbb{F}}^{f}$ is the canonical inclusion given by the embedding  
\begin{align*}
\mathbb{F}^{n}&\hookrightarrow \mathbb{F}^{\infty}\\
(x_{1},...,x_{n})&\mapsto (x_{1},...,x_{n},0,...,0).
\end{align*}  
Since $\Psi\circ\iota\circ\psi^{-1}$ is identity and $\psi$ is a bijection (see Theorem \ref{GeomeaningofVf}), the function $\Psi$ is onto. 

To see the injectivity of $\Psi$, we let $u=\{\mathfrak{U}_{\alpha},f_{\beta\alpha}\}$ and $v=\{\mathfrak{V}_{\gamma},g_{\delta\gamma}\}$ be two $iso\mathcal{V}_{\mathbb{F}}^{f}$-cocycles. Let $U$ and $V$ be the image of $u$ and $v$ under $\Psi$ and suppose they are isomorphic via the following bundle map
\begin{center}
\begin{tikzpicture}
\node(Lu) at (0,1) {$U$};
\node(Ru) at (2,1) {$V$};
\node(Ml) at (1,0) {$X$};

\path[->, font=\scriptsize,>=angle 90]

(Lu) edge node [above]{$T$}(Ru)  
(Lu) edge (Ml)  
(Ru) edge (Ml);

\end{tikzpicture}
\end{center}
Then $T$ induces a vector bundle isomorphism 
\[T_{\gamma\alpha}:f^{\ast}_{\alpha}\operatorname{E}_{k}\vert_{\mathfrak{U}_{\alpha}\cap\mathfrak{V}_{\gamma}}\rightarrow g^{\ast}_{\gamma}\operatorname{E}_{k}\vert_{\mathfrak{V}_{\gamma}\cap\mathfrak{U}_{\alpha}},\]
for each $\alpha$ and $\gamma$. For the sake of simplicity, we let $U_{\alpha,\gamma}$ denote $f^{\ast}_{\alpha}\operatorname{E}_{k}\vert_{\mathfrak{U}_{\alpha}\cap\mathfrak{V}_{\gamma}}$ and $V_{\gamma,\alpha}$ denote $g^{\ast}_{\gamma}\operatorname{E}_{k}\vert_{\mathfrak{V}_{\gamma}\cap\mathfrak{U}_{\alpha}}$.

Now, by Lemma \ref{Universalbundleofisom}, we know $T_{\gamma\alpha}$ induces a section 
\[t_{\gamma\alpha}:\mathfrak{U}_{\alpha}\cap\mathfrak{V}_{\gamma}\rightarrow iso_{k}(\mathbb{F})\] 
such that the following diagram commutes
\begin{center}
\begin{tikzpicture}
\node(Lu) at (0,2) {$\operatorname{Isom}(U_{\alpha,\gamma},V_{\gamma,\alpha})$};
\node(Ll) at (0,0) {$\mathfrak{U}_{\alpha}\cap\mathfrak{U}_{\gamma}$}; 
\node(Ru) at (4,2) {$iso_{k}(\mathbb{F})$};
\node(Rl) at (4,0) {$\operatorname{Gr}_{k}(\mathbb{F}^{\infty})\times\operatorname{Gr}_{k}(\mathbb{F}^{\infty})$};

\path[->, font=\scriptsize,>=angle 90]

(Lu) edge (Ru)  
(Lu) edge (Ll)
(Ll) edge node [above]{$(f_{\alpha},g_{\gamma})$}(Rl) 
(Ru) edge (Rl);

\draw [->](Ll) to [out=60,in=-60] node[right]{\scriptsize $T_{\gamma\alpha}$} (Lu);
\draw [->,dashed] (Ll) to node [above]{\scriptsize $t_{\gamma\alpha}$} (Ru);

\end{tikzpicture}
\end{center}
To see the collection $\{t_{\gamma\alpha}\}$ constitutes an isomorphism between $\{f_{\beta\alpha}\}$ and $\{g_{\delta\gamma}\}$, we need to show:
\begin{align*}
t_{\gamma\beta}\circ f_{\beta\alpha}&=t_{\gamma\alpha}\text{  on  } \mathfrak{U_{\alpha}\cap\mathfrak{U}_{\beta}\cap\mathfrak{V}_{\gamma}}\\
g_{\delta\gamma}\circ t_{\gamma\beta}&=t_{\delta\beta} \text{  on  }
\mathfrak{U_{\beta}\cap\mathfrak{V}_{\gamma}\cap\mathfrak{V}_{\delta}}
\end{align*} 
  
Now, observe the commutative diagram below  
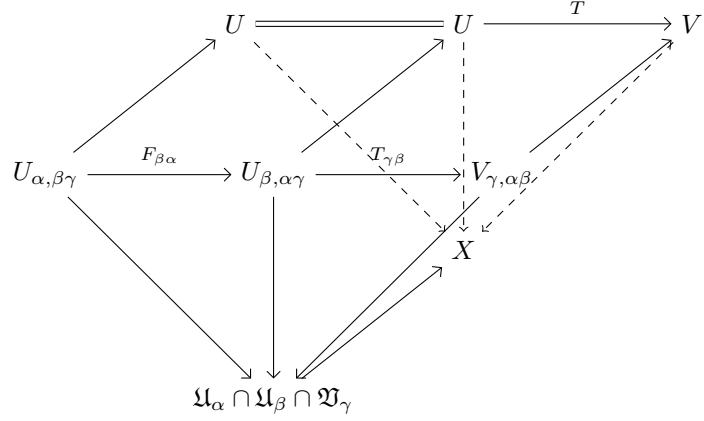
\begin{center}
\begin{tikzpicture}
\node(FLu) at (0,2) {$U_{\alpha,\beta\gamma}$};
\node(FMu) at (3,2) {$U_{\beta,\alpha\gamma}$};
\node(FRu) at (6,2) {$V_{\gamma,\alpha\beta}$};
\node(FMl) at (3,-1) {$\mathfrak{U}_{\alpha}\cap\mathfrak{U}_{\beta}\cap\mathfrak{V}_{\gamma}$};

\node(BLu) at (2.5,4) {$U$};
\node(BMu) at (5.5,4) {$U$};
\node(BRu) at (8.5,4) {$V$};
\node(BMl) at (5.5,1) {$X$};

\path[->, font=\scriptsize,>=angle 90]

(FLu) edge node [above]{$F_{\beta\alpha}$}(FMu)  
 
(FMu) edge node [above]{$T_{\gamma\beta}$}(FRu)
(FMu) edge (FMl)
(FRu) edge (FMl) 
(FLu) edge (FMl)

(FLu) edge (BLu)
(FMu) edge (BMu)
(FRu) edge (BRu)
(FMl) edge (BMl)

(BMu) edge node [above]{$T$}(BRu);
\draw [->,dashed] (BLu) to (BMl); 
\draw [->,dashed] (BMu) to (BMl);
\draw [->,dashed] (BRu) to (BMl);
 
\draw [double equal sign distance] (BLu) to (BMu); 
\end{tikzpicture}
\end{center}
where $U_{\alpha,\beta\gamma}$ (resp. $U_{\beta,\alpha\gamma}$ and $V_{\gamma,\alpha\beta}$) stands for the restriction of $U_{\alpha,\beta}$ (resp. $U_{\beta,\gamma}$ and $V_{\gamma,\beta}$) on $\mathfrak{U}_{\alpha}\cap\mathfrak{U}_{\beta}\cap\mathfrak{V}_{\gamma}$. The diagram tells us that $T_{\gamma\alpha}= T_{\gamma\beta}\circ F_{\beta\alpha}$, or equivalently, the front square (with bold arrows) in the diagram below commutes  
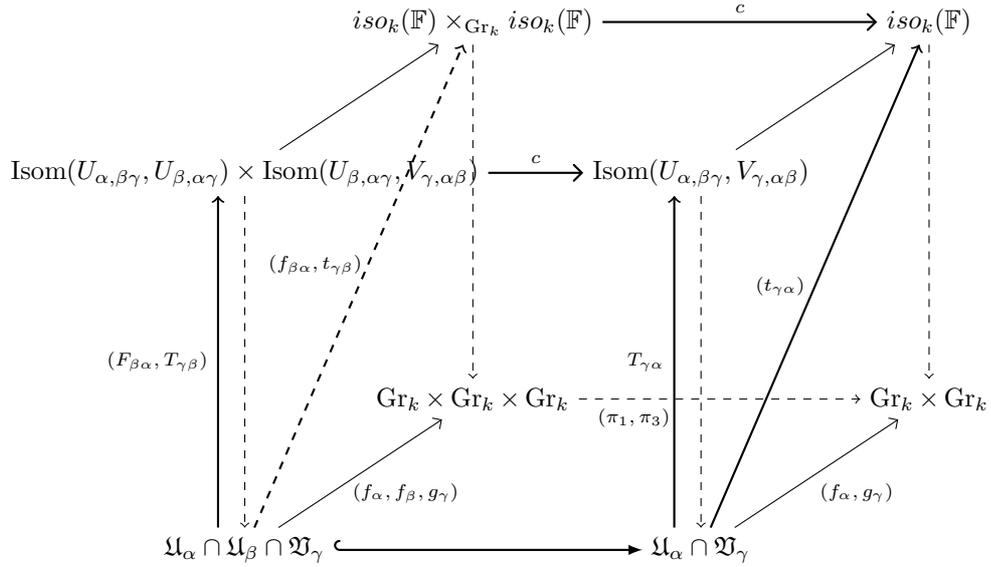
\begin{center}
\begin{tikzpicture}
\node(FLu) at (0,5) {$\operatorname{Isom}(U_{\alpha,\beta\gamma},U_{\beta,\alpha\gamma})\times \operatorname{Isom}(U_{\beta,\alpha\gamma},V_{\gamma,\alpha\beta})$};
\node(FLl) at (0,0) {$\mathfrak{U}_{\alpha}\cap\mathfrak{U}_{\beta}\cap\mathfrak{V}_{\gamma}$};

\node(FRu) at (6,5) {$\operatorname{Isom}(U_{\alpha,\beta\gamma},V_{\gamma,\alpha\beta})$};
\node(FRl) at (6,0) {$\mathfrak{U}_{\alpha}\cap\mathfrak{V}_{\gamma}$};

\node(BLu) at (3,7) {$iso_{k}(\mathbb{F})\times_{\operatorname{Gr}_{k}}iso_{k}(\mathbb{F})$};
\node(BLl) at (3,2) {$\operatorname{Gr}_{k}\times\operatorname{Gr}_{k}\times\operatorname{Gr}_{k}$};

\node(BRu) at (9,7) {$iso_{k}(\mathbb{F})$};
\node(BRl) at (9,2) {$\operatorname{Gr}_{k}\times\operatorname{Gr}_{k}$};

\path[thick,->, font=\scriptsize,>=angle 90]

(FLu) edge node [above]{$c$}(FRu)  

(BLu) edge node [above]{$c$}(BRu);  

\path[->, font=\scriptsize,>=angle 90]
(FLu) edge (BLu)
(FLl) edge node [below,xshift=1.7em]{$(f_{\alpha},f_{\beta},g_{\gamma})$}(BLl)
(FRu) edge (BRu)
(FRl) edge node [below,xshift=1.5em]{ $(f_{\alpha},g_{\gamma})$}(BRl);

\draw [->,dashed] (FRu) to (FRl); 
\draw [->,dashed] (FLu) to (FLl);
\draw [->,dashed] (BRu) to (BRl);
\draw [->,dashed] (BLu) to (BLl);
\draw [->,dashed] (BLl) to node [below,xshift=-3.2em]{\scriptsize $(\pi_{1},\pi_{3})$}(BRl); 
 
\draw [thick, right hook-latex](FLl) to (FRl); 

\draw [transform canvas={xshift=-1em}, thick, ->](FLl) to node[left]{\scriptsize $(F_{\beta\alpha},T_{\gamma\beta})$}(FLu);  
\draw [thick,dashed, ->](FLl) to node [left,yshift=.8em,xshift=.5em]{\scriptsize $(f_{\beta\alpha},t_{\gamma\beta})$}(BLu);
\draw [thick, ->](FRl) to node [left]{\scriptsize $(t_{\gamma\alpha})$} (BRu); 
\draw [transform canvas={xshift=-1em}, thick,->] (FRl) to node[left]{\scriptsize $T_{\gamma\alpha}$} (FRu);

\end{tikzpicture}
\end{center}
Hence, the square (bold (dashed) arrows) that slices the cube also commutes. In particular, we have proved
\[t_{\gamma\beta}\circ f_{\beta\alpha}=t_{\gamma\alpha}\text{  on  } \mathfrak{U_{\alpha}\cap\mathfrak{U}_{\beta}\cap\mathfrak{V}_{\gamma}}\]
The identity 
\[g_{\delta\gamma}\circ t_{\gamma\beta}=t_{\delta\beta}\text{  on  } \mathfrak{U_{\alpha}\cap\mathfrak{U}_{\beta}\cap\mathfrak{V}_{\gamma}}\] 
can be proved in a similar manner. Thus, the function $\Psi$ is bijective.
\end{proof}
 
\section{The relation between $\coprod\limits_{k}\operatorname{Gr}_{k}(\mathbb{F}^{\infty})$, $\mathcal{V}_{\mathbb{F}}$ and $\mathcal{V}_{\mathbb{F}}^{f}$} 
Let $\mathcal{G}$ be the category internal in $\mathpzc{Top}_{\ast}^{w}$ given by  
\[O(\mathcal{G})=M(\mathcal{G})=\coprod_{k}\operatorname{Gr}_{k}(\mathbb{F}^{\infty})\]
and trivial structure maps---namely, the identity. Then we have an obvious inclusion 
\[\mathcal{G}\hookrightarrow iso\mathcal{V}_{\mathbb{F}}^{f}.\]
Combing with the canonical inclusion
\[iso\mathcal{V}_{\mathbb{F}}\rightarrow iso\mathcal{V}_{\mathbb{F}}^{f}\]
given by the map
\begin{align*}
\mathbb{F}^{n}&\hookrightarrow \mathbb{F}^{\infty}\\
(x_{1},...,x_{n})&\mapsto (x_{1},...,x_{n},0,...,0),
\end{align*}
we have the following theorem.   
\begin{theorem}\label{ComparisonVFVFfG}
The cospan of inclusions
\[iso\mathcal{V}_{\mathbb{F}}\rightarrow iso\mathcal{V}_{\mathbb{F}}^{f}\leftarrow \mathcal{G}\]
induces a cospan of homotopy equivalences of $H$-semiring spaces
\begin{equation}\label{Cospanofheqofsemiring}
\vert \operatorname{Ner}_{\cdot}iso\mathcal{V}_{\mathbb{F}}\vert \rightarrow \vert \operatorname{Ner}_{\cdot} iso\mathcal{V}_{\mathbb{F}}^{f}\vert\leftarrow \vert \operatorname{Ner}_{\cdot}\mathcal{G}\vert=\coprod_{k}\operatorname{Gr}_{k}(\mathbb{F}^{\infty}). 
\end{equation}
\end{theorem}
\begin{proof} 
By the definitions, it is clear that the maps in the cospan \eqref{Cospanofheqofsemiring} preserves the additive and multiplicative identities. To see the map
\[\vert \operatorname{Ner}_{\cdot}iso\mathcal{V}_{\mathbb{F}}\vert \rightarrow \vert \operatorname{Ner}_{\cdot} iso\mathcal{V}_{\mathbb{F}}^{f}\vert\]
respects the additive and multiplicative operations, we observe that the diagrams 
\begin{center}
\begin{tikzpicture}
\node(Lu) at (0,2) {$iso\mathcal{V}_{\mathbb{F}}^{f}\times iso\mathcal{V}_{\mathbb{F}}^{f}$};
\node(Ll) at (0,0) {$iso\mathcal{V}_{\mathbb{F}}\times iso\mathcal{V}_{\mathbb{F}}$}; 
\node(Ru) at (7,2) {$iso\mathcal{V}_{\mathbb{F}}^{f}$};
\node(Rl) at (7,0) {$iso\mathcal{V}_{\mathbb{F}}$};

\path[->, font=\scriptsize,>=angle 90]

(Lu) edge node [above]{$\oplus $}(Ru)  
(Ll) edge (Lu)
(Ll) edge node [above]{$\oplus $}(Rl) 
(Rl) edge (Ru);

\end{tikzpicture}
\end{center} 
and
\begin{center}
\begin{tikzpicture}
\node(Lu) at (0,2) {$iso\mathcal{V}_{\mathbb{F}}^{f}\times iso\mathcal{V}_{\mathbb{F}}^{f}$};
\node(Ll) at (0,0) {$iso\mathcal{V}_{\mathbb{F}}\times iso\mathcal{V}_{\mathbb{F}}$}; 
\node(Ru) at (7,2) {$iso\mathcal{V}_{\mathbb{F}}^{f}$};
\node(Rl) at (7,0) {$iso\mathcal{V}_{\mathbb{F}}$};

\path[->, font=\scriptsize,>=angle 90]

(Lu) edge node [above]{$\otimes$}(Ru)  
(Ll) edge (Lu)
(Ll) edge node [above]{$\otimes$}(Rl) 
(Rl) edge (Ru);

\end{tikzpicture}
\end{center}
commute up to natural transformations. The natural transformations are induced by the isomorphisms  
\[\theta_{\ast}(\mathbb{F}^{k}\oplus \mathbb{F}^{l})\xleftarrow{\sim} \mathbb{F}^{k}\oplus \mathbb{F}^{l}=\mathbb{F}^{k}\oplus \mathbb{F}^{l}\xrightarrow[\eqref{DirectsumVf}]{\sim} \mathbb{F}^{k+l}\hookrightarrow \mathbb{F}^{\infty}\] 
and 
\[ \kappa_{\ast}(\mathbb{F}^{k}\otimes \mathbb{F}^{l}) \xleftarrow{\sim} \mathbb{F}^{k}\otimes \mathbb{F}^{l}=\mathbb{F}^{k}\otimes \mathbb{F}^{l}\xrightarrow[\eqref{TensorVf}]{\sim} \mathbb{F}^{kl},\]
respectively. Similarly, the diagram
\begin{center}
\begin{tikzpicture}
\node(Lu) at (0,2) {$iso\mathcal{V}_{\mathbb{F}}^{f}\times iso\mathcal{V}_{\mathbb{F}}^{f}$};
\node(Ll) at (0,0) {$\mathcal{G}\times\mathcal{G}$}; 
\node(Ru) at (7,2) {$iso\mathcal{V}_{\mathbb{F}}^{f}$};
\node(Rl) at (7,0) {$\mathcal{G}$};

\path[->, font=\scriptsize,>=angle 90]

(Lu) edge node [above]{$\oplus/\otimes$}(Ru)  
(Ll) edge (Lu)
(Ll) edge node [above]{$\oplus/\otimes$}(Rl) 
(Rl) edge (Ru);

\end{tikzpicture}
\end{center}
commutes up to natural transformations and hence the map
\[\coprod_{k}\operatorname{Gr}_{k}(\mathbb{F}^{\infty})\rightarrow \vert \operatorname{Ner}_{\cdot}iso\mathcal{V}_{\mathbb{F}}^{f}\vert\]
also respects the $H$-semiring space structures.

Now, given a paracompact Hausdorff space $X$, we observe the following diagram  
\begin{center}
\begin{tikzpicture}
\node(MT) at (3.5,3.5) {$\coprod\limits_{k}\operatorname{Vect}_{k}(X)$};
\node(Lu) at (0,2) {$[X,Biso\mathcal{V}_{\mathbb{F}}]$};
\node(Ll) at (0,0) {$[X,\vert\operatorname{Ner}_{\cdot}iso\mathcal{V}_{\mathbb{F}}\vert]$};
\node(Mu) at (3.5,2) {$[X,Biso\mathcal{V}_{\mathbb{F}}^{f}]$};
\node(Ml) at (3.5,0) {$[X,\vert\operatorname{Ner}_{\cdot}iso\mathcal{V}_{\mathbb{F}}^{f}\vert]$};
\node(Ru) at (7,2) {$[X,B\mathcal{G}]$};
\node(Rl) at (7,0) {$[X,\vert\operatorname{Ner}_{\cdot}\mathcal{G}\vert]$};

\path[->, font=\scriptsize,>=angle 90] 

(Ll) edge (Ml)
(Rl) edge (Ml)
(Lu) edge (Mu)
(Ru) edge (Mu)  
(Lu) edge (Ll)
(Mu) edge (Ml)  
(Ru) edge (Rl);
\draw[<->](MT) to node [above]{\scriptsize $b_{1}$}(Lu);
\draw[<->](MT) to node [right]{\scriptsize $b_{2}$}(Mu);
\draw[<->](MT) to node [above]{\scriptsize $b_{3}$}(Ru);

\end{tikzpicture}
\end{center} 
where the bijections $b_{1}$, $b_{2}$ and $b_{3}$ are given by Theorem $3.2$ in \cite{Wang4}. The vertical arrows are also bijective in view of Theorem $4.3$ in \cite{Wang4} and the properness of the simplicial spaces $\operatorname{Ner}_{\cdot}iso\mathcal{V}_{\mathbb{F}}, \operatorname{Ner}_{\cdot}iso\mathcal{V}_{\mathbb{F}}^{f}$ and $\operatorname{Ner}_{\cdot}\mathcal{G}$. In fact, the space $B\mathcal{G}$ is homeomorphic to $\vert\operatorname{Ner}_{\cdot}\mathcal{G}\vert\times \triangle^{\infty}$ and the map $\vert\operatorname{Ner}_{\cdot}\mathcal{G}\vert\times \triangle^{\infty}\simeq B\mathcal{G}\rightarrow \vert\operatorname{Ner}_{\cdot}\mathcal{G}\vert$ is the obvious projection. As a result, all horizontal maps are bijective as well, and this, the maps in the cospan \eqref{Cospanofheqofsemiring} are homotopy equivalences.
\end{proof}

For the last result, we need to recall the notion of prespectra \cite{MMSS}. A prespectrum is a sequence of topological spaces and maps $\mathbf{E}=\{E_{i},\sigma_{i}:E_{i}\wedge S^{1}\rightarrow E_{i+1}\}_{i\in\mathbb{N}\cup\{0\}}.$ When the adjoint of $\sigma_{i}$ is a homotopy equivalence, it is called an $\Omega$-prespectrum. Some properties of the category of prespectra and its comparison with the category of $\operatorname{CW}$-spectra are discussed in \cite[Appendix $A$]{Wang1}. Now, observe that the singular functor 
\[\operatorname{Sing}_{\cdot}:\mathpzc{Top}^{w}\rightarrow s\operatorname{Sets}\] 
preserves finite limit, so by Lemma \ref{Functorpreserveslimits}, it sends an internal category in $\mathpzc{Top}^{w}$ to an internal category in $s\operatorname{Sets}$, where $s\operatorname{Sets}$ is the category of simplicial sets. Let $\mathcal{V}_{\mathbb{F},\cdot}^{f}:=\operatorname{Sing}_{\cdot}\mathcal{V}_{\mathbb{F}}^{f}$. Then we have the following theorem.  
\begin{theorem}
The internal category $\mathcal{V}_{\mathbb{F},\cdot}^{f}$ is a simplicial $S$-category equipped with a sum functor \cite[Definition $1.8$]{Wang5}, where the sum functor is given by the functor
\[\oplus:\mathcal{V}_{\mathbb{F}}^{f}\times \mathcal{V}_{\mathbb{F}}^{f}\rightarrow \mathcal{V}_{\mathbb{F}}^{f}.\]  
Furthermore, its associated $\Omega$-prespectrum $\mathbf{K}^{S}(\mathcal{V}_{\mathbb{F},\cdot}^{f})$ represents connective topological $K$-theory in the sense of \cite[Definition $3.2$]{Wang5}. 
\end{theorem}
\begin{proof}
Applying the functor $\operatorname{Sing}_{\cdot}$ to Theorem \ref{HsemiringVFf}, we get the first assertion. As for the second assertion, we recall that a $(-1)$-connective $\Omega$-prespectrum $\mathbf{E}$ represents connective topological $K$-theory if and only if it is equipped with a map of prespectra 
\[\mathbf{B}:\mathbf{E}\rightarrow \Omega^{d}\mathbf{E},\]
and for any pointed compact Hausdorff space $X$, there is an isomorphism  
\[[X,E_{0}]_{\ast}\cong \tilde{K}^{\mathbb{F}}(X)\]
such that the following commutes
\begin{center}
\begin{equation}\label{CTKdiagram}
\begin{tikzpicture}[baseline=(current bounding box.center)]
\node(Lu) at (0,2) {$[X,E_{0}]_{\ast}$};
\node(Ll) at (0,0) {$\tilde{K}^{\mathbb{F}}(X)$}; 
\node(Ru) at (4,2) {$[X,\Omega^{d}E_{0}]_{\ast}$};
\node(Rl) at (4,0) {$\tilde{K}^{\mathbb{F}}(X\wedge S^{d})$};

\path[->, font=\scriptsize,>=angle 90] 

(Lu) edge node [above]{$B_{0,\ast}$}(Ru)  
(Lu) edge node [right]{$\wr$}(Ll)
(Ll) edge node [above]{$B_{\ast}$}(Rl) 
(Ru) edge node [right]{$\wr$}(Rl);
\end{tikzpicture}
\end{equation}
\end{center}
where $B_{\ast}$ is induced from the map \eqref{Bottmapspacelevel}. $\mathbf{E}$ is a connective complex (resp. real) topological $K$-theory prespectrum when $d=2$ (resp. $d=8$).

Now, as shown in the proof of Lemma \ref{Manyrealizations}, the simplicial space $\operatorname{Ner}_{\cdot}\mathcal{V}_{\mathbb{F}}^{f}$ is proper, and for each $i$, the space $\operatorname{Ner}_{i}\mathcal{V}_{\mathbb{F}}^{f}$ has the homotopy type of a $\operatorname{CW}$-complex, so the canonical map
\[\vert\operatorname{Ner}_{\cdot}\mathcal{V}_{\mathbb{F},\cdot}^{f}\vert\rightarrow \vert\operatorname{Ner}_{\cdot}\mathcal{V}_{\mathbb{F}}^{f}\vert\]
is a homotopy equivalence, and thus, there is a homotopy equivalence 
\begin{equation}\label{HeqbetwenKVFfandZG}
K^{S}(\mathcal{V}_{\mathbb{F},\cdot}^{f})\simeq \mathbb{Z}\times G(\mathbb{F}^{\infty})
\end{equation} 
by the uniqueness of the topological group completion of $\vert iso\mathcal{V}_{\mathbb{F}}^{f}\vert$ and the second homotopy equivalence of $H$-spaces in Theorem \ref{ComparisonVFVFfG}.

On the other hand, from the proof of Theorem \ref{HsemiringVFf}, we know the tensor product $\otimes$ induces a simplicial bi-exact functor \cite[Definition $1.21$]{Wang5}
\[\mathcal{V}_{\mathbb{F},\cdot}^{f}\times\mathcal{V}_{\mathbb{F},\cdot}^{f}\rightarrow \mathcal{V}_{\mathbb{F},\cdot}^{f}\]
and hence a pairing 
\[ K^{S}(\mathcal{V}^{f}_{\mathbb{F},\cdot})_{0}\wedge K^{S}(\mathcal{V}^{f}_{\mathbb{F},\cdot})_{0}\rightarrow K^{S}(\mathcal{V}^{f}_{\mathbb{F},\cdot})_{0}.\]
This pairing realizes the multiplicative structure in topological $K$-theory by Theorem $1.22$ and Lemma $3.1.iii$ in \cite{Wang5}; therefore, one can simply define the Bott map to be the adjoint of the composition
\[\mathbf{K}^{S}(\mathcal{V}_{\mathbb{F},\cdot}^{f})\wedge S^{d}\xrightarrow{\operatorname{id}\wedge b} \mathbf{K}^{S}(\mathcal{V}_{\mathbb{F},\cdot}^{f})\wedge  \mathbb{Z}\times G(\mathbb{F}^{\infty}) \xrightarrow{\operatorname{id}\wedge \eqref{HeqbetwenKVFfandZG}} \mathbf{K}^{S}(\mathcal{V}_{\mathbb{F},\cdot}^{f})\wedge  K^{S}(\mathcal{V}_{\mathbb{F},\cdot}^{f})_{0}\rightarrow \mathbf{K}^{S}(\mathcal{V}_{\mathbb{F},\cdot}^{f}),\]
where $b$ is the map representing the Bott element. The commutativity of the diagram \eqref{CTKdiagram} then follows automatically.
\end{proof}

\addtocontents{toc}{\protect\setcounter{tocdepth}{0}}
\begin{appendices}
\addtocontents{toc}{\protect\setcounter{tocdepth}{1}}
\section{Weak Hausdorff $k$-spaces}
Here we review some properties of the category of weak Hausdorff $k$-spaces $\mathpzc{Top}^{w}$. Most of them are taken either directly or indirectly from \cite{Du}, \cite{FP}, \cite{JEA}, \cite{Kie}, \cite{Lew}, \cite{Str1},\cite{Str2} and \cite{MP}. Explanations are given for those statements that are lass obvious.  

\subsection{Basic Definitions} 
We begin with a sequence of basic definitions. Let $\mathcal{T}$ denote the category of topological spaces. Then a space $X\in\mathcal{T}$ is a $k$-space if and only if every compactly open subset is open, where a subset $U\subset X$ is called compactly open if, for any continuous map $f:K\rightarrow X$ with $K$ compact, the preimage $f^{-1}(U)$ is open in $K$. A space $X$ is weak Hausdorff if and only if, for any map $f:K\rightarrow X$ with $K$ compact, the image $f(K)$ is closed in $X$. The category of weak Hausdorff $k$-spaces\footnote{It is sometimes called compactly generated spaces, for example, in \cite{Lew} and \cite{MP}.} is denoted by $\mathpzc{Top}^{w}$. The category $\mathpzc{Top}^{w}$ is large enough for our works; it includes the spaces encountered in this paper.  
\begin{example}\label{ExamplesofkandCG}
\begin{enumerate}
\item Every $\operatorname{CW}$-complex is in $\mathpzc{Top}^{w}$.
\item Every topological manifold\footnote{By definition, topological manifolds are Hausdorff.} is in $\mathpzc{Top}^{w}$. 
\item An open or closed subset of a weak Hausdorff $k$-space in $\mathcal{T}$ is a weak Hausdorff $k$-space.
\end{enumerate}
\end{example}
Though the colimits and limits in the categories $\mathcal{T}$ and $\mathpzc{Top}^{w}$ are in general different, in some cases, they do coincide---to distinguish, we use $\operatorname{lim}^{t}$ and $\operatorname{colim}^{t}$ to denote limits and colimits in $\mathcal{T}$, respectively, whereas $\operatorname{lim}$ and $\operatorname{colim}$ are reserved for limits and colimits in $\mathpzc{Top}^{w}$, respectively.

\begin{lemma}\label{Directlimitofinjections}
Given a direct system of weak Hausdorff $k$-spaces $\{X_{i},f_{i\leq j}:X_{i}\rightarrow X_{j}\}_{i,j\in I}$ with $f_{i\leq j}$ injective, for every $i\leq j\in I$, then 
\[\operatorname*{colim}_{i\in I}X_{i}=\operatorname*{colim^{t}}_{i\in I}X_{i},\]
or equivalently the colimit  
$\operatorname*{colim^{t}}\limits_{i\in I}X_{i}$
is a weak Hausdorff space.
\end{lemma}
\begin{proof}
See \cite[Proposition 9.3.b]{Lew}.
\end{proof} 

One important property of the category $\mathpzc{Top}^{w}$, which is not shared by the category $\mathcal{T}$, is that $\mathpzc{Top}^{w}$ is cartesian closed, namely that we have the adjoint functors
\[-\times X :\mathpzc{Top}^{w}\leftrightarrows \mathpzc{Top}^{w}:(-)^{X},\]
where $Y\times X$ is the product of $X$ and $Y$ in $\mathpzc{Top}^{w}$ and $Z^{X}$ is the function space from $X$ to $Z$ given by the compact-open topology and $k$-ification. In particular, the colimit in $\mathpzc{Top}^{w}$ commutes with finite product.

\begin{definition}
A (closed) subset $A\rightarrow X$ in $\mathpzc{Top}^{w}$ (resp. $\mathcal{T}$) is a (closed) cofibration\footnote{Any cofibration in $\mathpzc{Top}^{w}$ is closed.} if and only if, given any $Y\in\mathpzc{Top}^{w}$ and a commutative diagram below (solid arrows)  
\begin{center}
\begin{tikzpicture}
\node(Lu) at (0,2) {$A$};
\node(Ll) at (0,0) {$X$};
\node(Ru) at (2,2) {$Y^{I}$};
\node(Rl) at (2,0) {$Y$};

\path[->, font=\scriptsize,>=angle 90]

(Lu) edge (Ru) 
(Lu) edge (Ll)
(Ll) edge (Rl) 
(Ru) edge (Rl);
\draw[->,dashed] (Ll) to (Ru);
\end{tikzpicture}
\end{center}
the dashed arrow exists.

\noindent
A map $Y\rightarrow B$ in $\mathpzc{Top}^{w}$ (resp. $\mathcal{T}$) is a Hurewicz fibration if and only if, given any $X\in \mathpzc{Top}^{w}$ (resp. $\mathcal{T}$) and a commutative diagram below (solid arrows) 
\begin{center}
\begin{tikzpicture}
\node(Lu) at (0,2) {$X$};
\node(Ll) at (0,0) {$X\times I$};
\node(Ru) at (2,2) {$Y$};
\node(Rl) at (2,0) {$B$};

\path[->, font=\scriptsize,>=angle 90]

(Lu) edge (Ru) 
(Lu) edge (Ll)
(Ll) edge (Rl) 
(Ru) edge (Rl);
\draw[->,dashed] (Ll) to (Ru);
\end{tikzpicture}
\end{center}
the dashed arrow exists.
\end{definition}
\begin{theorem}
The category $\mathpzc{Top}^{w}$ is a proper monoidal model category whose weak equivalences, cofibrations and fibrations are cofibrations, Hurewicz fibrations and homotopy equivalences. 
\end{theorem}
\begin{proof}
\cite[Theorem 17.11]{MP}.
\end{proof}

\subsection{Pullbacks, colimits and cofibrations}
Here we recollect some special properties of the category $\mathpzc{Top}^{w}$. For convenience's sake, we denote the limit of a cospan $X\rightarrow Y\leftarrow Z$ in $\mathpzc{Top}^{w}$ by $X\times_{Y}Z$, and given a closed subspace $X_{0}\subset X$, $X\cup X_{0}\times I$ denotes the colimit of the span $X_{0}\times I\leftarrow X_{0}\rightarrow X$ in $\mathpzc{Top}^{w}$. 

\begin{lemma}\label{Colimitandpullback}
Let $I$ be a directed set. Given a cospan in the category of functors from $I$ to $\mathpzc{Top}^{w}$ $(X_{i},f_{i\leq j}:X_{i}\rightarrow X_{j})\rightarrow (Z_{i},h_{i\leq j}:Z_{i}\rightarrow Z_{j})\leftarrow (Y_{i},g_{i
\leq j}:Y_{i}\rightarrow Y_{j})$ with $Z_{i}\in \mathpzc{Top}^{w}$ and $f_{i\leq j},g_{i\leq j}$ and $h_{i\leq j}$ injective, for every $i\leq j\in I$, then we have the homeomorphism
\[\operatorname*{colim}_{i\in I}(X_{i}\times_{Z_{i}}Y_{i})=\operatorname*{colim}_{i\in I}X_{i}\times_{\operatorname*{colim}\limits_{i\in I}Z_{i}}\operatorname*{colim}_{i\in I}Y_{i} \]
\end{lemma}
\begin{proof} 
In view of the assumptions and Lemma \ref{Directlimitofinjections}, we can define  
\begin{align*}
X&:=\operatorname*{colim}_{i\in I}X_{i}=\operatorname*{colim^{t}}_{i\in I}X_{i};\\
Z&:=\operatorname*{colim}_{i\in I}Z_{i}=\operatorname*{colim^{t}}_{i\in I}Z_{i};\\
Y&:=\operatorname*{colim}_{i\in I}Y_{i}=\operatorname*{colim^{t}}_{i\in I}Y_{i},
\end{align*}
and on the other hand, there are homeomorphisms 
\[\operatorname*{colim}_{j\in I}\operatorname*{colim}_{i\in I}(X_{i}\times_{Z}Y_{j})=\operatorname*{colim}_{i\in I}(X_{i}\times_{Z}Y_{i})=\operatorname*{colim}_{i\in I}(X_{i}\times_{Z_{i}}Y_{i}).\]
Then, applying \cite[Corollary $10.4$]{Lew} to the present case, we get the following:
\begin{multline*}  
\operatorname*{colim}_{i\in I}(X_{i}\times_{Z_{i}}Y_{i})=\operatorname*{colim}_{j\in I}\operatorname*{colim}_{i\in I}(X_{i}\times_{Z}Y_{j})=\operatorname*{colim^{t}}_{j\in I}\operatorname*{colim^{t}}_{i\in I}(X_{i}\times_{Z}Y_{j})\\
\overset{\cite{Lew}}{=}\operatorname*{colim^{t}}_{j\in I}(X\times_{Z} Y_{j})=X\times_{Z}Y. 
\end{multline*} 
\end{proof}
 
\begin{lemma}\label{PullbackandCof}
\begin{itemize}
\item A map $A\rightarrow X$ is a cofibration in $\mathpzc{Top}^{w}$ if and only if it is closed cofibration in $\mathcal{T}$.

\item Given a pullback diagram in $\mathpzc{Top}^{w}$
\begin{center}
\begin{tikzpicture}
\node (Lu) at (0,1){$A$};
\node (Ll) at (0,0){$B$};
\node (Rl) at (2,0){$Y$};
\node (Ru) at (2,1){$X$};

\path[->, font=\scriptsize,>=angle 90] 

(Lu) edge node[above] {$j$}(Ru)
(Lu) edge node[right] {$q$}(Ll)
(Ru) edge node[right] {$p$}(Rl)
(Ll) edge node[above] {$i$}(Rl);
\end{tikzpicture}
\end{center}
with $i$ a cofibration and $p$ a Hurewicz fibration in $\mathpzc{Top}^{w}$, then the induced maps $j$ and $q$ are cofibration and Hurewicz fibration in $\mathpzc{Top}^{w}$, respectively.
\end{itemize}
\end{lemma}
\begin{proof}
A map is a closed cofibration $i:A\rightarrow X$ in $\mathcal{T}$ if and only if there is a retract $X\times I\rightarrow X\times\{0\}\cup A\times I$ \cite[Theorem 2]{Str1}, but closed cofibrations in $\mathpzc{Top}^{w}$ implies the latter and thus it is also a closed cofibration in $\mathcal{T}$. 

In view of the first assertion, we know $A\rightarrow X$ is a cofibration in $\mathpzc{Top}^{w}$ if and only if there exists a map $\phi:X\rightarrow I$ with $\phi^{-1}(0)=A$ and a homotopy $H:X\times I\rightarrow X$ such that $H(x,0)=\operatorname{id}$, $H(a,t)=a$, for every $a\in A$, and $H(x,t)\in A$ when $t>\phi(x)$ \cite[Lemma $4$]{Str1}. With this observation, the second statement can be proved in exactly the same way as in \cite[Theorem 12]{Str1}.  
\end{proof}

The next lemma is Kieboom's theorem \cite[Theorem $1$]{Kie} in the category of weak Hausdorff $k$-spaces.
\begin{lemma}\label{PullbackandCof2}
Consider the commutative diagram in $\mathpzc{Top}^{w}$
\begin{center}
\begin{tikzpicture}
\node(Lu) at (0,1) {$X_{0}$};
\node(Ll) at (0,0) {$X$};
\node(Mu) at (1.5,1) {$B_{0}$};
\node(Ml) at (1.5,0) {$B$};
\node(Ru) at (3,1) {$Y_{0}$};
\node(Rl) at (3,0) {$Y$};

\path[->, font=\scriptsize,>=angle 90]

(Ll) edge (Ml)
(Rl) edge node [above]{$p$}(Ml)
(Lu) edge (Mu)
(Ru) edge node [above]{$p_{0}$}(Mu)  
(Lu) edge node [right]{$i$}(Ll)
(Mu) edge node [right]{$j$}(Ml)  
(Ru) edge node [right]{$k$}(Rl);
\end{tikzpicture}
\end{center}
in which $p_{0}$ and $p$ are Hurewicz fibrations and $i$, $j$ and $k$ are cofibrations. Then the map
\[X_{0}\times_{B_{0}} Y_{0}\rightarrow X\times_{B} Y\]
is also a cofibration in $\mathpzc{Top}^{w}$.
\end{lemma}
\begin{proof}
Kieboom's pullback theorem asserts the same statement in the category $\mathcal{T}$, but his method can be carried over to the category $\mathpzc{Top}^{w}$. Following Kieboom's proof \cite[p381-382]{Kie}, we decompose the map
\[X_{0}\times_{B_{0}} Y_{0}\rightarrow X\times_{B} Y\]
into two maps
\begin{align}
X_{0}\times_{B_{0}} Y_{0}&\rightarrow X_{0}\times_{B_{0}} Y\vert_{B_{0}};\label{PTCof1}\\
X_{0}\times_{B_{0}} Y\vert_{B_{0}}&\rightarrow X\times_{B} Y,\nonumber 
\end{align}
where $Y\vert_{B_{0}}$ is the pullback of $Y\rightarrow B$ along $B_{0}$. By Lemma \ref{PullbackandCof}, the second map is a cofibration because the Hurewicz fibration $X_{0}\times_{B_{0}}Y\vert_{B_{0}}\rightarrow X_{0}$ is homeomorphic to the pullback of the Hurewicz fibration $X\times_{B} Y\rightarrow  X$ along $X_{0}\rightarrow X$. To see the first map is also a closed cofibration, we first note that Lemma \ref{PullbackandCof} implies $Y\vert_{B_{0}}\rightarrow Y$ is a cofibration and hence $Y_{0}\rightarrow Y\vert_{B_{0}}$ is also a cofibration \cite[Lemma $5$]{Str2}---the lemma in \cite{Str2} is stated in the category $\mathcal{T}$, yet in view of Lemma \ref{PullbackandCof}, the same assertion holds in the category $\mathpzc{Top}^{w}$. Secondly, by Lemma $3.8$ in \cite{JEA} , we know there is a retraction 
\begin{equation}\label{RetractionLA7}
r:E\vert_{B_{0}}\times I\rightarrow Y\vert_{B_{0}}\cup Y_{0}\times I
\end{equation} 
over $B_{0}$ \cite[Definition $1.1$]{HK}. Again, in \cite{HK} and \cite{JEA}, they are working in the category $\mathcal{T}$. However, their approaches can be applied also to the category $\mathpzc{Top}^{w}$. In more details, one observes that, similar to the case in $\mathcal{T}$ \cite[p.325]{JEA}, a map $E\xrightarrow{q} L$ in $\mathpzc{Top}^{w}$ is a Hurewicz fibration in $\mathpzc{Top}^{w}$ if and only if there exists a universal section of the map 
\[E^{I}\xrightarrow{\gamma\mapsto (\gamma(0),q\circ \gamma)} E\times_{L} L^{I},\]    
where $E\times_{L} L^{I}$ is given by the limit 
\[\operatorname{lim}(E\xrightarrow{q}L   \xleftarrow{\operatorname{ev}_{0}}L^{I})\]
and $\operatorname{ev}_{0}(\beta)=\beta(0)$. This observation implies that one can attach Hurewicz fibrations in $\mathpzc{Top}^{w}$ as in $\mathcal{T}$ \cite[Lemma $3.8$]{JEA}---the argument in \cite[Lemma]{JEA} works also in the category $\mathpzc{Top}^{w}$. Therefore, the map $Y\vert_{B_{0}}\cup Y_{0}\times I\rightarrow B_{0}$ is also a Hurewicz fibration in $\mathpzc{Top}^{w}$, and the diagram below
\begin{center}
\begin{tikzpicture}
\node (Lu) at (0,1){$Y\vert_{B_{0}}\cup Y_{0}\times I$};
\node (Ll) at (0,0){$Y\times I$};
\node (Rl) at (4,0){$B_{0}$};
\node (Ru) at (4,1){$Y\vert_{B_{0}}\cup Y_{0}\times I$};

\path[->, font=\scriptsize,>=angle 90] 

(Lu) edge (Ru)
(Lu) edge (Ll)
(Ru) edge (Rl)
(Ll) edge (Rl);
\draw[->,dashed](Ll) to node [above] {\scriptsize $r$} (Ru);
\end{tikzpicture}
\end{center}
gives us the retraction \eqref{RetractionLA7}.

Applying Proposition $7.4$ and Proposition $10.4$, the retraction $r$ further induces a retraction 
\[X_{0}\times_{B_{0}}E\vert_{B_{0}}\times I\rightarrow X_{0}\times_{B_{0}}E\vert_{B_{0}}\cup X_{0}\times_{B_{0}}E_{0}\times I\] 
which implies the map \eqref{PTCof1} is a closed cofibration.
 
\end{proof}

\begin{lemma}\label{Colimitandcof}
Given a levelwise cofibration between two sequences of closed cofibrations in $\mathpzc{Top}^{w}$
\[\{X_{i},f_{i}\}\xrightarrow{\{h_{i}\}}\{Y_{i},g_{i}\},\]
we have the colimit
\[\operatorname*{colim}_{i}h_{i}:\operatorname*{colim}_{i}X_{i}\rightarrow \operatorname*{colim}_{i}Y_{i}\]
is also a closed cofibration.
\end{lemma}
\begin{proof}
Recall that the category of functors from $\mathbb{N}$ to $\mathpzc{Top}^{w}$, denoted by $\mathpzc{Top}^{w,\mathbb{N}}$, can be endowed with the Reedy model structure, in which a map    
\[A=\{A_{i},a_{i}\}\xrightarrow{\{c_{i}\}}B=\{B_{i},b_{i}\},\]
is a cofibration if and only if the map from the latching object $L_{i}B\rightarrow B_{i}$ is a closed cofibration \cite[p.45]{Du}, \cite[Section VII.2]{GJ} for every $i$. Now, by our assumption and Lillig's union theorem \cite[p.410]{Li}, we know the map of sequences
\[\{X_{i},f_{i}\}\xrightarrow{\{h_{i}\}}\{Y_{i},g_{i}\},\] 
constitutes a cofibration in this model structure, and hence, the lemma follows from the Quillin adjunction \cite[Section $8.8$]{Du}
\[\operatorname{colim}:\mathpzc{Top}^{w,\mathbb{N}}\leftrightarrows \mathpzc{Top}^{w}:c.\]
\end{proof}
\subsection{Colimit and subspace topologies}
Given a sequence of spaces $\{Y_{i},h_{i}\}$ in $\mathpzc{Top}^{w}$ and a subspace $X\subset Y:=\operatorname*{colim}\limits_{i}Y_{i}$ in $\mathpzc{Top}^{w}$, the subspace topology on $X$ is in general different from the colimit topology induced from $X_{i}:= X\cap Y_{i}$. In the following situation, these two topologies coincide, however.
\begin{lemma}\label{Colimandsubspace}
Let $\{Y_{i},f_{i}:Y_{i}\rightarrow Y_{i+1}\}$ be an object in $\mathpzc{Top}^{w,\mathbb{N}}$ with $Y_{i}$ normal and $f_{i}$ a cofibration, for every $i\in\mathbb{N}$. Now, given a subspace $X$ of $\operatorname*{colim}\limits_{i\in\mathbb{N}}Y_{i}$ such that the subspace $X_{i}:=X\cap Y_{i}\subset X$ in $\mathcal{T}$ is a $k$-space. Then we have 
\[\operatorname*{colim}_{i\in\mathbb{N}}X_{i}=X\in\mathpzc{Top}^{w},\]
namely that the subspace topology of $X$ in $\mathcal{T}$ is identical to the colimit topology of $X$ induced by $X_{i}$. Note that, a priori, we do not know if $X\in\mathpzc{Top}^{w}$.
\end{lemma} 
\begin{proof}
The assertion follows from \cite[Proposition A.5.4]{FP} where it states, as objects in $\mathcal{T}$, we have the homeomorphism
\[\operatorname*{colim^{t}}_{i\in\mathbb{N}}X_{i}=X.\] However, by the assumption $X_{i}\in\mathpzc{Top}^{w}$ and Lemma \ref{Directlimitofinjections}, we know this homeomorphism actually holds in $\mathpzc{Top}^{w}$. 
\end{proof}

\subsection{Hurewicz fibrations and spaces of the homotopy type of $\operatorname{CW}$-complexes}
A detailed explanation of the following proposition can be found in \cite[Chapter $5$]{FP}
\begin{lemma}\label{TCWfibration}
\begin{enumerate} 
\item Given a Hurewicz fibration $X\rightarrow Y\rightarrow Z$ in $\mathpzc{Top}^{w}$, if $X$ and $Z$ have the homotopy type of $\operatorname{CW}$-complexes, then so is $Y$.
\item Given a Hurewicz fibration $X\rightarrow Y\rightarrow Z$ in $\mathpzc{Top}^{w}$, if $Y$ and $Z$ have the homotopy type of $\operatorname{CW}$-complexes, then so is $X$.
\end{enumerate}
\end{lemma}   
\begin{proof}
See \cite[Proposition 5.4.1, Theorem 5.4.2]{FP}.
\end{proof} 
In general, given a Hurewicz fibration $X\rightarrow Y\rightarrow Z$ in $\mathpzc{Top}^{w}$ with $X$ and $Y$ having the homotopy type of $\operatorname{CW}$-complexes, $Z$ does not necessarily have the homotopy type of a $\operatorname{CW}$-complex.

\section{Internal categories}
In this section, we review some basic properties of categories internal in a category with finite limits.    

\begin{definition}[Category internal in $\mathcal{M}$]\label{ICT}
Let $\mathcal{M}$ be a category with finite limits. Then a category internal in $\mathcal{M}$, denoted by $\mathcal{C}$, is a collection consisting of two objects $O(\mathcal{C})$ and $M(\mathcal{C})$ and four morphisms    
\begin{align*}
s_{c}: &M(\mathcal{C})\rightarrow O(\mathcal{C});\\
t_{c}: &M(\mathcal{C})\rightarrow O(\mathcal{C});\\
e_{c}: &O(\mathcal{C})\rightarrow M(\mathcal{C});\\
\circ_{c}:& M(\mathcal{C})\times_{O(\mathcal{C})}M(\mathcal{C})\rightarrow M(\mathcal{C}),
\end{align*}
in $\mathcal{M}$ such that they satisfy the following commutative diagrams
\begin{center} 
\begin{tikzpicture}
\node (Lu) at(0,2){$O(\mathcal{C})$}; 
\node (Ru) at(2,2){$M(\mathcal{C})$};
\node (Rl) at(2,0){$O(\mathcal{C})$};

\path[->, font=\scriptsize,>=angle 90]
 
(Lu) edge node[above]{$e_{c}$}(Ru)
(Ru) edge node[right]{$s_{c}/ t_{c}$}(Rl)
(Lu) edge node[left]{$\operatorname{id}$}(Rl);

\node (Lu) at(5,2){$M(\mathcal{C})\times_{O(\mathcal{C})}M(\mathcal{C})$}; 
\node (Ll) at(5,0){$M(\mathcal{C})$};
\node (Ru) at(8,2){$M(\mathcal{C})$};
\node (Rl) at(8,0){$O(\mathcal{C})$};

\path[->, font=\scriptsize,>=angle 90]
(Lu) edge node [right]{$p_{1}/p_{2}$} (Ll)
(Lu) edge node [above]{$\circ_{c}$}(Ru)
(Ll) edge node [below]{$s_{c}/ t_{c}$}(Rl)
(Ru) edge node [right]{$s_{c}/ t_{c}$}(Rl);

\end{tikzpicture}
\end{center}
\begin{center}
\begin{tikzpicture}
\node (Lu) at(1,2){$M(\mathcal{C})\times_{O(\mathcal{C})}M(\mathcal{C})\times_{O(\mathcal{C})}M(\mathcal{C})$}; 
\node (Ll) at(1,0){$M(\mathcal{C})\times_{O(\mathcal{C})}M(\mathcal{C})$};
\node (Ru) at(7,2){$M(\mathcal{C})\times_{O(\mathcal{C})}M(\mathcal{C})$};
\node (Rl) at(7,0){$M(\mathcal{C})$};

\path[->, font=\scriptsize,>=angle 90]
(Lu) edge node [left]{$\circ_{c}\times_{O(\mathcal{C})}\operatorname{id}$} (Ll)
(Lu) edge node [above]{$\operatorname{id}\times_{O(\mathcal{C})}\circ_{c}$}(Ru)
(Ll) edge node [below]{$\circ_{c}$}(Rl)
(Ru) edge node [right]{$\circ_{c}$}(Rl);

\node (Lu) at(0,-1.5){$O(\mathcal{C})\times_{O(\mathcal{C})}M(\mathcal{C})$}; 
\node (Ru) at(8,-1.5){$M(\mathcal{C})\times_{O(\mathcal{C})}O(\mathcal{C})$};
\node (Mu) at(4,-1.5){$M(\mathcal{C})\times_{O(\mathcal{C})}M(\mathcal{C})$};
\node (Ml) at(4,-3.5){$M(\mathcal{C})$};

\path[->, font=\scriptsize,>=angle 90]
(Lu) edge node [above]{$e_{c}\times_{O(\mathcal{C})}\operatorname{id}$} (Mu)
(Ru) edge node [above]{$\operatorname{id}\times_{O(\mathcal{C})}e_{c}$}(Mu)
(Mu) edge node [right]{$\circ_{c}$}(Ml)
(Ru) edge node [right,yshift=-.1cm]{$p_{1}$}(Ml)
(Lu) edge node [left,yshift=-.1cm]{$p_{2}$}(Ml);
\end{tikzpicture}
\end{center}
where $p_{1}$ and $p_{2}$ are obvious projections to the first component and the second component, respectively. The objects $O(\mathcal{M})$ and $M(\mathcal{M})$ are called objects of objects and morphisms, respectively; the morphisms
$s_{c}$, $t_{c}$, $e_{c}$ and $\circ_{c}$ are called the source morphism, target morphism, identity-assigning morphism and the composition of $\mathcal{C}$, respectively.    
\end{definition}
In this paper, we are only concerned with the case $\mathcal{M}=\mathpzc{Top}^{w}$ or $s\operatorname{Sets}$.

\begin{definition}[Internal Functor]\label{IF}
A(n) (internal) functor of categories internal in $\mathcal{M}$  
\[\mathcal{F}:\mathcal{C}\rightarrow\mathcal{D},\]
consists of two morphisms $O(\mathcal{F}),M(\mathcal{F})\in M(\mathcal{M})$ such that they are compatible with the structure maps, namely the commutative diagrams
\begin{center}
\begin{tikzpicture}
\node (Lu) at(0,2){$M(\mathcal{C})$}; 
\node (Ll) at(0,0){$O(\mathcal{C})$};
\node (Ru) at(2,2){$M(\mathcal{D})$};
\node (Rl) at(2,0){$O(\mathcal{D})$};

\path[->, font=\scriptsize,>=angle 90]
(Lu) edge node [left]{$s_{c}/t_{c}$} (Ll)
(Lu) edge node [above]{$M(\mathcal{F})$}(Ru)
(Ll) edge node [below]{$O(\mathcal{F})$}(Rl)
(Ru) edge node [right]{$s_{d}/t_{d}$}(Rl);

\node (Lu) at(4,2){$O(\mathcal{C})$}; 
\node (Ll) at(4,0){$M(\mathcal{C})$};
\node (Ru) at(6,2){$O(\mathcal{D})$};
\node (Rl) at(6,0){$M(\mathcal{D})$};

\path[->, font=\scriptsize,>=angle 90]
(Lu) edge node [left]{$e_{c}$} (Ll)
(Lu) edge node [above]{$O(\mathcal{F})$}(Ru)
(Ll) edge node [below]{$M(\mathcal{F})$}(Rl)
(Ru) edge node [right]{$e_{d}$}(Rl);
\end{tikzpicture}
\end{center}

\begin{center}
\begin{tikzpicture}
\node(Lu) at (0,2) {$M(\mathcal{C})\times_{O(\mathcal{C})}M(\mathcal{C})$};
\node(Ll) at (0,0) {$M(\mathcal{D})\times_{O(\mathcal{D})}M(\mathcal{D})$}; 
\node(Ru) at (4,2) {$M(\mathcal{C})$};
\node(Rl) at (4,0) {$M(\mathcal{D})$};

\path[->, font=\scriptsize,>=angle 90]

(Lu) edge node [above]{$\circ_{c}$}(Ru)  
(Lu) edge node [right]{$M(\mathcal{F})\times_{O(\mathcal{F})}M(\mathcal{F})$}(Ll)
(Ll) edge node [above]{$\circ_{d}$}(Rl) 
(Ru) edge node [right]{$M(\mathcal{F})$}(Rl);

\end{tikzpicture}

\end{center} 
\end{definition}

\begin{lemma}\label{Functorpreserveslimits}
Given a functor $\mathcal{F}:\mathcal{M}\rightarrow \mathcal{N}$ 
that preserves finite limits, then $\mathcal{F}$ induces a functor from the category of internal categories in $\mathcal{M}$ to the category of internal categories in $\mathcal{N}$.
\end{lemma}  
\begin{proof}
Observe that the commutative diagrams in \ref{ICT} and \ref{IF} involve only finite limits.
\end{proof}
The main example in this paper is the singular functor $\operatorname{Sing}_{\cdot}:\mathpzc{Top}^{w}\rightarrow s\operatorname{Sets}$.

\begin{definition}[Natural Transformation]\label{NT}
A(n) (internal) natural transformation of two internal functors of categories internal in $\mathcal{M}$, 
\[\mathcal{F};\mathcal{G}:\mathcal{C}\rightarrow\mathcal{D},\] 
is a morphism $\phi:O(\mathcal{C}) \rightarrow M(\mathcal{D})\in M(\mathcal{M})$ such that the following diagrams commute:
\begin{center}
\begin{tikzpicture}
\node (Lu) at(0,2){$O(\mathcal{C})$}; 
\node (Ru) at(2,2){$M(\mathcal{D})$};
\node (Rl) at(2,0){$O(\mathcal{D})$};

\path[->, font=\scriptsize,>=angle 90]
(Lu) edge node [above]{$\phi$}(Ru)
(Lu) edge node [left]{$O(\mathcal{F})/O(\mathcal{G})$}(Rl)
(Ru) edge node [right]{$s_{d}/t_{d}$}(Rl);

\node (Mu) at(7,2){$M(\mathcal{C})$}; 
\node (Lm) at(5,1){$M(\mathcal{D})\times_{O(\mathcal{D})}M(\mathcal{C})$};
\node (Rm) at(9,1){$M(\mathcal{D})\times_{O(\mathcal{D})}M(\mathcal{D})$};
\node (Ml) at(7,0){$M(\mathcal{D})$};

\path[->, font=\scriptsize,>=angle 90]
(Mu) edge node [above]{$\hat{G}$} (Lm)
(Mu) edge node [above]{$\hat{F}$}(Rm)
(Lm) edge node [below,yshift=-.1em]{$\circ_{c}$}(Ml)
(Rm) edge node [below,yshift=-.15em,xshift=.1em]{$\circ_{d}$}(Ml);

\end{tikzpicture}
\end{center}
where the morphisms $\hat{F}$ and $\hat{G}$ are induced from the following diagrams.
\begin{center} 
\begin{tikzpicture}
\node (Mu) at(2,2){$M(\mathcal{C})$}; 
\node (Lm) at(0,1){$M(\mathcal{D})$};
\node (Rm) at(4,1){$M(\mathcal{D})$};
\node (Ml) at(2,0){$O(\mathcal{D})$};

\path[->, font=\scriptsize,>=angle 90]
(Mu) edge node [left, yshift=.1cm]{$\phi\circ s_{c}$} (Lm)
(Mu) edge node [right,yshift=.1cm]{$M(\mathcal{G})$}(Rm)
(Lm) edge node [below, yshift=-.1em]{$t_{d}$}(Ml)
(Rm) edge node [below, yshift=-.15em,xshift=.1em]{$s_{d}$}(Ml);

\node (Mu) at(8,2){$M(\mathcal{C})$}; 
\node (Lm) at(6,1){$M(\mathcal{D})$};
\node (Rm) at(10,1){$M(\mathcal{D})$};
\node (Ml) at(8,0){$O(\mathcal{D})$};

\path[->, font=\scriptsize,>=angle 90]
(Mu) edge node [left,yshift=.1cm]{$\phi\circ t_{c}$} (Lm)
(Mu) edge node [right,yshift=.1cm]{$M(\mathcal{F})$}(Rm)
(Lm) edge node [below,yshift=-.1em]{$s_{d}$}(Ml)
(Rm) edge node [below,yshift=-.1em,xshift=.1em]{$t_{d}$}(Ml);

\end{tikzpicture}
\end{center}
\end{definition}

If $\mathcal{M}=\mathpzc{Top}^{w}$, then there is an alternative description of a(n) (internal) natural transformation.
 
\begin{lemma}\label{EqDefforNaTrans}
Let $[1]$ be the category of $\lbrace 0,1\rbrace$ with morphisms $\lbrace 0\rightarrow 0,0\rightarrow 1, 1\rightarrow 1\rbrace$ and consider the product internal category 
\[\mathcal{C}\times [1],\] 
whose space of objects is 
\[O(\mathcal{C})\times\{0\} \coprod O(\mathcal{C})\times\{1\}\] 
and space of morphisms
\[M(\mathcal{C})\times \{0\rightarrow 0\}\coprod M(\mathcal{C})\times\{0\rightarrow 1\}\coprod M(\mathcal{C})\times\{1\rightarrow 1\},\]
then given a(n) (internal) natural transformation (Definition \ref{NT}) is equivalent to given a(n) (internal) functor 
\[\Phi:\mathcal{C}\times [1]\rightarrow \mathcal{D}\] 
with 
\begin{align*}
\Phi\vert_{O(\mathcal{C})\times\{0\}}&=   O(\mathcal{F})\\
\Phi\vert_{O(\mathcal{C})\times\{1\}}&=   O(\mathcal{G})\\
\Phi\vert_{M(\mathcal{C})\times\{0\rightarrow 0\}}& = M(\mathcal{F})\\
\Phi\vert_{M(\mathcal{C})\times\{1\rightarrow 1\}}& =  M(\mathcal{G})
\end{align*}
\end{lemma}
\begin{proof}
The second commutative diagram in Definition \ref{NT} is equivalent to the following commutative diagram
\begin{center}
\begin{tikzpicture}
\node(Lm) at (0,2) {$M(\mathcal{C})$};
\node(Ml) at (3,0)   {$M(\mathcal{D})\times_{O(\mathcal{D})}M(\mathcal{D})$}; 
\node(Mu) at (3,4)   {$M(\mathcal{D})\times_{O(\mathcal{D})}M(\mathcal{D})$};
\node(Rm) at (6,2) {$M(\mathcal{D})$};

\draw [->] (Lm) to[out=-90,in=180] node [right,yshift=.3cm]{\scriptsize $(\phi\circ s_{c}, M(\mathcal{G}))$}(Ml);
\draw [->] (Ml) to[out=0,in=-90] node [left]{\scriptsize $\circ_{d}$}(Rm); 
\draw [->] (Mu) to[out=0,in=90] node [left]{\scriptsize $\circ_{d}$}(Rm); 
\draw[->] (Lm) to [out=90,in=180]  node [right,yshift=-.3cm]{\scriptsize $(M(\mathcal{F}),\phi\circ t_{c})$}(Mu);
\end{tikzpicture}
\end{center}
Thus, we can define $\Phi$ as follows:
\[\Phi(h,0\rightarrow 1):=\mathcal{F}(h)\circ \phi(t_{d}(f))= \phi(s_{c}(f))\circ \mathcal{G}(h).\]
For the converse, we let $\phi(c):=\Phi(\operatorname{id}_{c},0\rightarrow 1)$.
\end{proof}

In the following, we recall the nerve construction of an internal category in $\mathcal{M}$, which is a functor from the category of internal categories in $\mathcal{M}$ to the category of simplicial objects in $\mathcal{M}$. 

\begin{construction}\label{ConstinIntCat}
\begin{enumerate}
\item
Given $\mathcal{C}$ a category internal in $\mathcal{M}$, the nerve construction is given by 
\[\operatorname{Ner}_{k}\mathcal{C}:= \overbrace{M(\mathcal{C})\times_{O(\mathcal{C})}M(\mathcal{C})\times_{O(\mathcal{C})}.....\times_{O(\mathcal{C})}M(\mathcal{C})}^{k \text{ copies of }M(\mathcal{C})},\]
with the degeneracy map $\mathbcal{s}_{i}$ given by
\[\operatorname{id}\times_{O(\mathcal{C})}\operatorname{id}\times_{O(\mathcal{C})}...\overbrace{\times_{O(\mathcal{C})}e_{c}\times_{O(\mathcal{C})}}^{i}...\times_{O(\mathcal{C})}\operatorname{id},\]
and the face map $d_{i}$ by 
\[\operatorname{id}\times_{O(\mathcal{C})}\operatorname{id}\times_{O(\mathcal{C})}...\overbrace{\times_{O(\mathcal{C})}\circ_{c}\times_{O(\mathcal{C})}}^{i}...\times_{O(\mathcal{C})}\operatorname{id},\] 
where $i\neq 0,k$. $d_{0}, d_{k}$ are given by forgetting the first and the last component respectively.

\item An internal functor of internal categories in $\mathcal{M}$, 
\[\mathcal{F}:\mathcal{C}\rightarrow\mathcal{D},\] 
induces a simplicial map 
\[\operatorname{Ner}_{\cdot}\mathcal{F}:\operatorname{Ner}_{\cdot}\mathcal{C}\rightarrow \operatorname{Ner}_{\cdot}\mathcal{D},\] 
where 
\[\operatorname{Ner}_{k}\mathcal{F}:=\overbrace{M(\mathcal{F})\times_{O(\mathcal{F})}M(\mathcal{F})\times_{O(\mathcal{F})}...\times_{\mathcal{O}}M(\mathcal{F})}^{k}.\]

\item An internal natural transformation of internal functors in $\mathpzc{Top}^{w}$
\[\Phi:\mathcal{C}\times [1]\rightarrow \mathcal{D},\]
induces a homotopy  
\[\operatorname{Ner}_{\cdot}\mathcal{C}\times \triangle^{1}_{\cdot}\rightarrow \operatorname{Ner}_{\cdot}\mathcal{D},\]
where $\triangle^{1}_{k}:=\triangle([k],[1])$, the set of morphisms from $[k]$ to $[1]$ in the simplex category $\triangle$.  
\end{enumerate}
\end{construction} 
  
\end{appendices}

\addcontentsline{toc}{section}{R\hspace*{1.8em}References}

\bibliographystyle{alpha}

\bibliography{Reference}

\end{document}